\documentclass[review]{elsarticle}
\usepackage{fullpage}
\usepackage{lineno,hyperref}
\modulolinenumbers[5]

\usepackage{mathtools}
\usepackage{amsfonts}

\usepackage{hyperref}
\usepackage{cleveref}
\usepackage{color}
\usepackage{graphicx}
\usepackage{subcaption}
\usepackage{url}
\usepackage{listings}

\usepackage{codehighlight}
\usepackage{pifont}

\usepackage{placeins}

\let\Oldsection\section
\renewcommand{\section}{\FloatBarrier\Oldsection}

\let\Oldsubsection\subsection
\renewcommand{\subsection}{\FloatBarrier\Oldsubsection}

\let\Oldsubsubsection\subsubsection
\renewcommand{\subsubsection}{\FloatBarrier\Oldsubsubsection}

\usepackage{svg}

\usepackage{todonotes}

\newcommand{\dx}{\text{d}x}

\newcommand{\dt}{\Delta t}
\newcommand{\norm}[1]{\lVert#1\rVert}

\begin{document}

\begin{frontmatter}

\title{Hybrid FEM-NN models: Combining artificial neural networks with the finite element method}

\author[simula]{Sebastian K. Mitusch\corref{mycorrespondingauthor}}
\cortext[mycorrespondingauthor]{Corresponding author}
\ead{sebastkm@simula.no}

\author[simula]{Simon W. Funke}
\ead{simon@simula.no}

\author[simula]{Miroslav Kuchta}
\ead{miroslav@simula.no}

\address[simula]{Simula Research Laboratory, 1364 Fornebu, Norway}

\begin{abstract}
We present a methodology combining neural networks with physical principle constraints 
in the form of partial differential equations (PDEs). The approach allows to train neural 
networks while respecting the PDEs as a strong constraint in the optimisation as apposed to making them part of the loss 
function. The resulting models are discretised in space by the finite element method (FEM). The method 
applies to both stationary and transient as well as linear/nonlinear PDEs. 
We describe implementation of the approach as an extension of the existing 
FEM framework FEniCS and its algorithmic differentiation tool dolfin-adjoint.
Through series of examples we demonstrate capabilities of the approach to recover coefficients and 
missing PDE operators from observations. Further, the proposed method is compared with alternative 
methodologies, namely, physics informed neural networks and standard PDE-constrained 
optimisation. Finally, we demonstrate the method on a complex cardiac cell model problem using deep neural networks.
\end{abstract}


\end{frontmatter}

    \section{Introduction}

    The plummeting cost of physical sensors, computational power, and data storage results
    in an explosion of data, and the task of manually extracting information from that
    data has become overwhelming. In the last decade, statistical learning, and
    specifically artificial neural networks (NN), have proven to be immensely valuable
    in meeting these challenges. One disadvantage of neural networks is however that
    a priori they embed no inherent knowledge of the physical, or mathematical, laws
    governing the underlying systems at hand \cite{marcus2018deep}. Embedding such
    knowledge is non trivial, requiring e.g. novel architectures \cite{long2019pde, long2018pde, ruthotto2019deep, lu2019deeponet, li2020fourier}
    or problem formulations \cite{raissi_physics-informed_2019}.

    Conversely, models based on physical principles, typically described by  partial/ordinary differential equations
    (PDEs/ODEs), have been employed for centuries; these models have the advantages of a solid mathematical
    foundation and a wide array of numerical methods facilitating their solution.
    However, constructing PDE models, and pursuant solution techniques, is an arduous task.
    PDE models are often rigid, relying on explicit assumptions, or are so large as to be
    computationally infeasible. The primary limitation of PDE models is a lack of ability
    to learn new principles from observational input.

    Bridging the gap between explicit PDE systems and observation-driven NN learning to overcome the limitations of either approaches has recently gained scientific traction.
    Raissi et al. \cite{raissi_physics-informed_2019} introduced physics-informed neural networks (PINNs) for solving PDEs by training a neural network with a loss function consisting of the PDE residual.
    They further demonstrate that PINNs can be combined with discovery of governing equations through identification of coefficients in the PDE.
    This equation discovery technique was later demonstrated with a neural network approximator \cite{tartakovsky2018learning, raissi_deep_2018}.
    PINNs are able to approximate solutions even of complex non-linear PDEs, which is essential when using a neural network approximator for equation discovery.
    
    However, when combining equation discovery with PINNs the training problem becomes more complex.
    Instead of only training to find an approximation for the unknown equation terms, the training process also involves approximating the solution to this equation.
    This additional training step might be unnecessary if traditional solvers for PDEs are used to solve the approximated equation.
    Furthermore, traditional PDE solvers such as the finite element method (FEM) or finite volume method have rich theory,
    providing convergence guarantees. In comparison, such theories for PINNs have only
    recently emerged e.g. \cite{shin2020convergence, wang2020and, mishra2020estimates}.

    {
    In this paper, we propose a methodology to combine the finite element method with neural networks.
    The aim is to obtain an approach that combines the strengths of PDE-based modelling, e.g.
    the geometric flexibility and rich set of finite element functions of the finite element method,
    with the flexibility of neural networks to express unknown functions.
    Equation discovery using traditional PDE solvers in conjunction with neural
    networks has previously been demonstrated.
    In \cite{berg_neural_2017} a neural network with a single hidden layer is
    used to approximate a spatially varying diffusion coefficient in the Poisson
    equation that was discretised with the finite element method.
    The framework ADCME.jl presented in \cite{huang_learning_2020, xu2020adcme}
    enables training of finite element/finite difference-discretized models
    with neural networks and has been succesfully applied e.g. to identify physical parameters
    in the steady state Navier-Stokes equation \cite{fan_solving_2020} or
    learning the constitutive relations \cite{huang_learning_2020, xu2020inverse}.
    If used with the finite elements ADCME.jl
    requires that the user implements backpropagation pass through the
    layer(s) representing the discretization and solution of the PDE constraint.
    In contrast, in our work the pass is automatically generated based on symbolic
    representation of variational form of the constraint. 
    Rackauckas et. al \cite{rackauckas_universal_2020} propose Universal Differential Equations in which ODEs and PDEs are augmented with neural networks on finite difference schemes.
    The results show promise for the development of hybrid PDE-NN models, but do not demonstrate the combination of more flexible PDE discretisation methods such as the FEM.
    }

    Here we present a novel approach for training neural networks
    augmenting PDEs, which addresses some of the shortcomings
    of the previous works. In particular, the approach applies to both
    stationary and transient as well as linear/non-linear PDEs.
    The hybrid PDE-NN models are discretised in space using FEM, enabling the use of the well-established finite element framework FEniCS \cite{logg_dolfin_2010}.
	The neural networks can either be defined directly in FEniCS or through the machine learning library PyTorch \cite{NEURIPS2019_9015}.
	We demonstrate the approach on a variety of problems, including problems with partial observations, noisy observations and deep neural networks.
	Our work is structured as follows. In section \ref{sec:framework}
    the hybrid methodology is introduced in a general setting along
    with details about discretisation and the resulting optimisation problem.
	Afterwards, section \ref{sec:examples} presents a range of examples of increasing complexity.
	Finally, section \ref{sec:conclusion} summarizes and provides discussion on future work.
    
        \section{{Framework for hybrid FEM-NN models}}\label{sec:framework}
        Let $\Omega\subset\mathbb{R}^n$, $n\geq 1$, $T>0$. In the following we
        consider a class of PDE-constrained optimization problems of the form
        \begin{equation}\label{eq:general}
        \min_{{u}, N}\sum_{i=1}^nL(u_i, {u}) + R({u}, N) \mbox{ subject to } E(u, N) = 0 \mbox{ in } \Omega\times(0, T),
        \end{equation}
        where $\mathcal{U} = \{u_i\}_{i=1}^{n}$ is a collection of observations of the state
        ${u}$ and $L$, $R$ are respectively the loss function and a regulatization. Note that only the observations of the state are assumed. 
        The operator $E$ is a differential operator encoding prior (partial) knowledge about
        the problem at hand in the form of a partial differential equation\footnote{Equipped with suitable boundary conditions.}
        to be satisfied by ${u}$ and the unknown term $N$. We remark
        that $N$ can simply be a missing coefficient of the equation or
        an unknown differential operator (e.g. closure of a turbulence
        model \cite{turbulence}). In the former case \cref{eq:general} is a standard problem, that is, both $u$ and $N$ can be approximated as e.g. finite element functions.
        In the latter case, the problem can be translated into standard setting
        e.g. by using sparse regression technique \cite{Brunton3932}. Then, however,
        successful identification of the structure of $N$ requires additional
        knowledge in the form of a suitable dictionary of building expressions.
        In order to address both cases simultaneously we choose here to represent
        $N$ as a neural network $\mathcal{N}(u; W)$ parameterized by weights
        $W$. This choice gives rise to a (training) problem
        \begin{equation}\label{eq:general_NN}
        \min_{{u}, W}\sum_{i=1}^nL(u_i, {u}) + R({u}, \mathcal{N}(u;W)) \mbox{ subject to } E(u, \mathcal{N}(u;W)) = 0 \mbox{ in } \Omega\times(0, T).
        \end{equation}

        Commonly when the physics is only partially known, the unknown
        physical component can be more easily formulated in an additive way.
        Additionally, we mainly focus on time-dependent problems, from which
        we can formulate a specific class of \emph{hybrid PDE-NN} models
        \begin{equation}\label{eq:general_coupled_form}
          \begin{aligned}
    	    u_t &= F(u) + \mathcal{N}(u; W) \quad &&  \text{ in } \Omega \times (0, T),\\
    	    u &= u_0 && \text{ on } \Omega \times \{t = 0\},\\
            u &= u_b && \text{ on } \partial\Omega \times (0, T)
          \end{aligned}
        \end{equation}
        with solution $u: \Omega \times (0, T) \rightarrow \mathbb{R}^d$ defined on the
        space-time cylinder $\Omega \times(0, T)$.
        Here, $F$ represents an ordinary or partial differential operator,
        and $\mathcal{N}(u; W):  \Omega \times (0, T)  \rightarrow \mathbb{R}^d$ an artificial neural network with weights $W$. For simplicity, we assume that the neural networks only use point evaluations of the solution $u$. That is the neural network takes the form $\mathcal{N}(u; W)(x, t) \equiv  \mathcal{N}(x,  t, u(x, t), \nabla u (x, t), \dots; W)$. 

        We close this section by giving some concrete examples of the hybrid
        PDE-NN models. To this end, let us consider a heat equation with a
        nonlinear diffusion coefficient $\kappa$, that is, $u:\mathbb{R}\to\mathbb{R}$ satisfies
        $u_t-\nabla\cdot(\kappa(u)\nabla u)=0$ in $\Omega\times(0, T)$
        along with suitable initial and boundary conditions. Based on
        observations of $u$ we may now consider the problem of recovering
        the coefficient; a possible hybrid PDE-NN model then reads
        \[
    	    0 = E(u, N) := u_t + G(u, \mathcal{N}) \quad \text{ in } \Omega \times (0, T),\\
        \]
        where $G(u, \mathcal{N}) = -\nabla\cdot(\mathcal{N}(u; W)\nabla u)$ and $\mathcal{N}(u; W):\Omega \times (0, T)\to\mathbb{R}$. Thus in this case, $\mathcal{N}$ 
        approximates the diffusion coefficient.
        Alternatively, we may have $F(u, \mathcal{N}) = -\nabla\cdot(\mathcal{N}(u; W))$ where $\mathcal{N}(u; W):\Omega \times (0, T)\to\mathbb{R}^d$
        represents the flux of $u$. We remark that in the latter case our
        prior knowledge of the structure of the problem is less as $E$ simply encodes a balance law.

	\subsection{Discretisation in time}\label{sec:discretisation-in-time}
	We consider two strategies for discretising \cref{eq:general_coupled_form} in time:
        a direct discretisation and operator splitting.
        Using for instance a $\theta$-method in the first approach yields for each timestep
        a problem
	\begin{equation}
	u^{n+1} - u^n = \Delta t \left(F(u^{n+\theta}) + \mathcal{N}(u^{n+\theta}; W)\right)\quad\mbox{ in }\Omega,
	\label{eq:general_coupled_discrete_form_implicit}
    \end{equation}
    where $\Delta t$ denotes the time step and $u^{n+\theta}:=(1-\theta)u^{n} + \theta u^{n+1}$.
    Depending on the choice of $\theta$ one obtains an explicit Euler ($\theta=0$), semi-implicit Crank-Nicolson ($\theta=1/2$) or implicit Euler ($\theta=1$) method.
   	This approach has some difficulties in practice.
        The first challenge is to solve the highly non-linear problem \cref{eq:general_coupled_discrete_form_implicit}
        that arises at each time-step for $\theta>0$. The non-linearity stems primarily from the
        activation functions in the neural network $\mathcal{N}$ and will increase with the
        complexity of the neural network architecture. As a result, the convergence radius of
        standard non-linear solvers such as Newton's method, can be very small and might force
        the usage of very small time-steps or advanced non-linear solvers. The second difficulty
        is to solve the linear subproblems that arise during each non-linear iteration. The dimension
        of these linear problems is typically large and tailored preconditioners for iterative linear solvers 
        are required to solve them efficiently. However, the left hand side of the linear problem
        contains both a linearised PDE operator and a linearised neural network. Preconditioners for
        such linear problems have not been developed yet.
        We apply the direct approach \cref{eq:general_coupled_discrete_form_implicit} to examples 
        presented further in sections \ref{sec:examples:comparison} and \ref{sec:examples:advection}.
    
        In the operator splitting approach \cref{eq:general_coupled_form} is decomposed
        into simpler problems that can be treated individually using specialized numerical algorithms.
        For instance, applying the Marchuk-Yanenko \cite{GLOWINSKI20033} splitting to \cref{eq:general_coupled_form} yields two subproblems for
        each timestep:
        \begin{equation}    	\label{eq:general_splitted_discrete_form_implicit_nn}
    \begin{aligned}
    	\bar{u}^{n+1} &= u^{n} + \Delta t F(\bar{u}^{n+\theta}) &&\mbox{in }\Omega, \\
    	u^{n+1} &= \bar u^{n+1} + \Delta t \mathcal{N}(\bar{u}^{n+1}; W) &&\mbox{in }\Omega
    \end{aligned}
    \end{equation}
    with $\bar{u}^{n+\theta}=\theta\bar{u}^{n+1} + (1-\theta)u^n$, $0\leq \theta\leq 1$.
    Note that the first subproblem corresponds to a pure PDE problem, while the second subproblem reduces to a residual neural network step. 	We remark that the Marchuk-Yanenko splitting scheme is only first order accurate.
	Thus, it can be preferable to use Strang splitting \cite{GLOWINSKI20033} which is formally second-order accurate for sufficiently smooth operators $F$ and $\mathcal{N}$.

	The operator splitting approach overcomes most of the difficulties that arise with the first strategy. 
	Firstly, the PDE and neural network operators are
	different in nature and hence require different numerical solution techniques. 
    For instance, a natural choice is to treat the PDE subproblems implictly, and the 
    NN subproblem explicitly, potentially even with substepping.
	Secondly, splitting enables the usage of specialised numerical methods to each subproblem. 
	This is particularly important for the first subproblem, for which there exist well-established discretisation schemes, such as the finite element method, as well as fast iterative solvers and preconditioners.
	In addition, the splitting approach allows us to use two specialized software frameworks of highly efficient and well-tested numerical methods, e.g. FEniCS \cite{logg_dolfin_2010} for the PDE part and 
	PyTorch \cite{NEURIPS2019_9015} for the NN part.
	The splitting approach \cref{eq:general_splitted_discrete_form_implicit_nn} is applied to 
	the problem studied in section \ref{sec:examples:cardiac}.

	\subsection{Discretisation in space}
 	We employ the finite element method for the spatial discretisation of the hybrid PDE-NN problems  \cref{eq:general_coupled_discrete_form_implicit} and
        \cref{eq:general_splitted_discrete_form_implicit_nn}. For this reason we term the resulting discrete system the \emph{hybrid FEM-NN models}. 
        
        Let $\Omega_h$ be a triangulation of the domain
        $\Omega$ with $h>0$ a characteristic size of the elements and let $U, V$ be suitable, discrete finite element
        spaces of trial and test functions. Then, the variational form of \cref{eq:general_coupled_discrete_form_implicit}
        reads: Find $u^{n+1}\in U$ such that for all $v\in V$
	\begin{align}
	\int_\Omega (u^{n+1} - u^n  )\cdot v \dx = \Delta t \int_\Omega  F(u^{n+\theta})\cdot v \dx  +  \Delta t \int_\Omega \mathcal{N}(u^{n+\theta}; W)\cdot v \dx.
	\label{eq:discretisation:space:coupled}
        \end{align}
        We remark that the finite element spaces and functions depend on the mesh size parameter,
        however, for simplicity of notation we do not write $u_h$, $v_h$ or $U_h$, $V_h$
        and only use the subscript when the dependence needs to be highlighted.

    For the Marchuk-Yanenko splitting scheme \cref{eq:general_splitted_discrete_form_implicit_nn}, the variational problem become: First, find $\bar{u}^{n+1}\in U$ such that
    \begin{equation}
    	\int_\Omega (\bar{u}^{n+1} - u^n  )\cdot v \dx = \Delta t \int_\Omega  F(\bar{u}^{n+\theta})\cdot v \dx
    \end{equation}
    for all $v\in V$. Then, find $u^{n+1}\in U$ such that 
    \begin{equation}
    	\int_\Omega (u^{n+1} - \bar u^{n+1}  )\cdot v \dx = \Delta t \int_\Omega  \mathcal{N}(\bar{u}^{n+1}; W)\cdot v \dx
    	\label{eq:discretisation:space:splitting}
    \end{equation}
    for all $v\in V$. If $\mathcal{N}(\bar{u}^{n}; W)\in U$ the second problem can be solved pointwise as
    \begin{equation}
    	u^{n+1}(x) = \bar u^{n+1}(x) + \Delta t \mathcal{N}(\bar{u}^{n}; W)(x) .
    \end{equation}
    
    Note that both \cref{eq:discretisation:space:coupled} and \cref{eq:discretisation:space:splitting} require integration of a neural network over the spatial domain $\Omega$.
    This can be a major challenge as there is currently no theory on optimal quadrature rules for NNs,
    and deriving the analytical expression for these integrals is practically infeasible as the network grows \cite{kharazmi_variational_2019}.
    Thus, in order to compute these integrals we use a Gaussian quadrature rule (on the triangulation $\Omega_h$) with a
    degree based on the finite element discretisation scheme used to represent the state and NN.
	\subsection{Training procedure and adjoint equation}
	While in regular NN problems one often has direct data on the NN output,
	in our methodology the direct NN output is assumed unknown and only
        data on what the solution of the FEM-NN model should match are provided.
	Thus, the cost function to be minimised measures the error of the NN output only indirectly through the solution of the FEM-NN model.

	Traditionally, neural networks are trained using back-propagation, which involves computing
	the gradient of the cost using reverse mode algorithmic differentiation (AD).
	In reverse mode AD the gradient is computed by the adjoint of an expansion through the chain rule \cite{naumann_art_2012},
	which in our problem setting would consist of computing $\frac{\partial \hat u_i}{\partial W}^*[s]$.
	Here $s$ is the result of previous gradient computations, typically $s = \frac{\mathrm{d} L}{\mathrm{d} \hat u}$ and $^{*}$ denotes the Hermitian 
	adjoint. 
	Under certain assumptions, this adjoint operator can be efficiently computed using the adjoint equation. 
	Efficiently in that the computational cost does not scale multiplicatively with the number of weights in the neural network.
        
	Let us assume that the operator $E(\hat u, W) = \hat u_t - F(\hat u) + N(\hat u; W)$ is continuously
        differentiable, with a unique solution $\hat u$ for each configuration of weights $W$.
	Thus, there is a solution operator $\hat u(W)$ such that $E(\hat u(W), W) = 0$.
	Further assuming that the linearised operator $\partial E(\hat u(W), W)/\partial \hat u$ is continuously invertible,
        by the implicit function theorem, the solution operator $\hat u(W)$ is continuously differentiable and
        its derivative is given through the following equation \cite{Ulbrich2009}
	\begin{align*}
		\frac{\partial E(\hat u(W), W)}{\partial \hat u} \frac{\partial \hat u(W) }{\partial W} + \frac{\partial E(\hat u(W), W)}{\partial W} = 0.
	\end{align*}
	It follows that the adjoint operator $\frac{\partial \hat u_i}{\partial W}^*[s]$ is given by
	\begin{align*}
		\frac{\partial \hat u_i}{\partial W}^*[s] = -\frac{\partial E(\hat u(W), W)}{\partial W}^* \frac{\partial E(\hat u(W), W)}{\partial \hat u}^{-*}[s].
	\end{align*}
	Introducing an auxiliary variable named the adjoint state $\lambda$, the gradients for the PDE solution are computed in two steps:
	First by solving the adjoint equation
	\begin{align*}
		\frac{\partial E(\hat u(W), W)}{\partial \hat u}^{*} \lambda_s = s.
	\end{align*}
	Second by using the acquired adjoint state to compute the gradient contributions from the PDE solution
	\begin{align*}
		\frac{\partial \hat u_i}{\partial W}^*[s] = -\frac{\partial E(\hat u(W), W)}{\partial W}^* \lambda_s.
	\end{align*}
	Note that while the equation $E(\hat u(W), W) = 0$ may be nonlinear, the adjoint equation is always linear.
	
	In the case of a nonlinear PDE $E(\hat u(W), W) = 0$ an iterative scheme such as Newton's method or Picard iterations is employed.
	This solution process involves solving many linearised systems, which (hopefully) approach the true solution of the nonlinear PDE.
	Assuming that the final solution is close to the true solution,
	the adjoint equation is solved only once for the nonlinear PDE.
	Thus, backpropagation through a nonlinear PDE will involve as many or fewer linear solves than the forward computations.
	
	With the adjoint method we have a way to efficiently compute gradients of the now \emph{unconstrained} optimisation problem
	\begin{align*}
		\min_W \sum_{i=1}^{n} L(\hat u_i(W), u_i).
	\end{align*}
	From the perspective of the back-propagation algorithm, the PDE solver is just as any other arithmetic operation which we can differentiate.
	This enables the integration with existing machine learning frameworks such as PyTorch \cite{NEURIPS2019_9015}.
	
	\subsection{Implementation}
	In order to solve the PDEs with FEM we employ the finite element framework FEniCS \cite{logg_dolfin_2010} which allows PDEs to be specified through the domain specific language UFL \cite{alna_es_unified_2014}.
	Using UFL the user defines the PDE through its discrete variational formulation with close to mathematical notation.
	For example, defining the variational formulation of the PDE $F := -\Delta u - f = 0$ amounts to:
        \vspace{-10pt}        
        \begin{python}
F = grad(u)*grad(v)*dx - f*v*dx
        \end{python}
        
	After defining the variational form, it can be used with the FEniCS interface to automatically assemble and solve the system.
	Using dolfin-adjoint \cite{mitusch_dolfin-adjoint_2019} the PDE solution can be automatically differentiated.
	Exploiting the symbolic representation of the variational formulation of the PDE, in addition to the symbolic differentiation capabilities of UFL, dolfin-adjoint derives the UFL expression for the adjoint equation and solves it using the same discretisation as the original PDE solver.
	
	The NNs can be defined directly in UFL, allowing the FEniCS backend to take care of the evaluation at integration points.
	For example, a single hidden layer of a NN $v_1 = \sigma(W_1 x + W_2)$ with weights $W_1$, bias $W_2$, and activation function
        $\sigma = \tanh$ can be expressed in UFL as
        \vspace{-10pt}        
	\begin{python}
sigma = ufl.tanh
v_1 = sigma(W_1 * ufl.as_vector([x]) + W_2)\end{python}
	However, UFL does not have built-in support for vectorized activation functions.
	Thus, if the output is not scalar, one has to loop over the output\footnote{
          We remark that this stage concerns building a symbolic representation
          of $\sigma$ and the code generated for the function's evaluation
          might employ vectorization and other optimizations \cite{homolya2017exposing}.
          }:
        \vspace{-10pt}
	\begin{python}
def sigma(v_0):
    return ufl.as_vector([ufl.tanh(v_0[i]) for i in range(v_0.ufl_shape[0])])
        \end{python}
        In a similar way one can easily build UFL expressions for arbitrary layered feed-forward NNs,
        using the class \texttt{ANN} in the snippet below. The \texttt{ANN} class can be found at
        \url{https://github.com/sebastkm/hybrid-fem-nn}.
        
        In turn, a hybrid FEM-NN training problem can be defined
        in a few lines of code. For example, the Poisson problem from section \ref{sec:examples:poisson}
        with an assumed unknown spatially varying source $f$ term approximated by a NN with 10 hidden neurons
        and scalar output reads:
\vspace{-10pt}
	\begin{python}
from fenics import *
from dolfin_adjoint import *

from neural_network import ANN

# Setup mesh, FE space V with boundary conditions `bcs`
# Let observations be `obs`
hat_u = Function(V)

layers = [2, 10, 1]
bias = [True, True]
x, y = SpatialCoordinate(mesh)
net = ANN(layers, bias=bias, mesh=mesh)
E = inner(grad(u), grad(v))*dx - net([x, y])*v*dx

# Solve PDE
solve(lhs(E) == rhs(E), hat_u, bcs)

# L^2 error as loss
loss = assemble((hat_u - obs)**2*dx)

# Define reduced formulation of problem
hat_loss = ReducedFunctional(loss, net.weights_ctrls())

# Use scipy L-BFGS optimiser
opt_theta = minimize(hat_loss, method="L-BFGS-B")\end{python}
 	This pure UFL approach allows for the whole implementation to only rely on FEniCS and dolfin-adjoint.
 	If we instead wish to approximate $f$ using a NN implemented outside of UFL, for example in PyTorch,
	we proceed by interpolating the NN to a finite element space. The details of this approach will be
        discussed in a future paper.

    \section{Examples}\label{sec:examples}
    In this section, we use numerical examples to evaluate the capabilities and limitations of hybrid FEM-NN models.
    First, we explore how the numerical errors of the hybrid model solution depends on factors such as finite element discretisation, noise or partially lacking observations (section \cref{sec:examples:poisson}).
    We then compare and benchmark the method against two alternative approaches: pointwise estimator and physics-informed neural networks
    (\cref{sec:examples:comparison}). 
    A key feature of the hybrid approach is that the neural network can   learn entire unknown terms in the PDE. This is demonstrated in  \cref{sec:examples:advection}, where a missing advection term in the advection-diffusion equation is recovered.
    Finally, in \cref{sec:examples:cardiac} the method is demonstrated on a more challenging problem from the field of cardiac electrophysiology.
    The source code is available at \url{https://github.com/sebastkm/hybrid-fem-nn-examples}.

   As training algorithm, we mainly rely on L-BFGS \cite{liu1989limited} with line-search. This choice was made because our problems have low amounts of data, do not require mini-batching and early testing revealed that L-BFGS achieves superior convergence compared to stochastic gradient descent.
   In addition, we sometimes supplement L-BFGS with the SciPy \cite{2020SciPy-NMeth} implementation of truncated Newton (TNC) \cite{nash1984newton}, as it occasionally surpassed the performance of only L-BFGS.
    
    We let $\hat u$ denote the numerical state solution, and $u$ the true state solution, on the spatio-temporal domain $\Omega \times \mathcal{T}$ with $\mathcal{T} = (0, T)$.
    In order to measure the accuracy of $\hat u$, we define the  (relative) error measures:
   	\begin{equation}
   	  e_{\Omega \times \mathcal{T}}^2(\hat{u}, u) = \frac{\int_\mathcal{T} \int_\Omega \lvert\hat u - u\rvert^2 \, \dx \, \mathrm{d}t }{\int_\mathcal{T} \int_\Omega \lvert u \rvert^2 \, \dx \, \mathrm{d}t}
          \quad\mbox{ and }\quad
          e_{\Omega}^2(\hat{u}, u, t) = \frac{\int_\Omega \lvert\hat u(\cdot, t) - u(\cdot, t)\rvert^2 \, \dx }{\int_\Omega \lvert u(\cdot, t) \rvert^2 \, \dx},          
   	\end{equation}
        where the latter local-in-time measure is used to give insight
        about the error in the spatial domain at a specific time point. We note that if the 
        arguments of the measures are clear from the context we shall simply write 
        $e_{\Omega \times \mathcal{T}}$ or $e_{\Omega(t)}$.
        Similarly, we define measures for the accuracy of the learned sub-physic term. Naturally, these metrics can only be used for synthetic examples where the true sub-physics term is known. We denote $G(u)$ the true unknown term that the neural network is representing and define relative sub-physics errors:
   	\begin{equation}
   	  \mathcal{E}^2(u) = \frac{\int_\mathcal{T} \int_\Omega \lvert \mathcal{N}(u; W) - G(u) \rvert^2 \, \dx \, \mathrm{d}t }{\int_\mathcal{T} \int_\Omega \lvert G(u) \rvert^2 \, \dx \, \mathrm{d}t}
          \quad\mbox{ and }\quad
          \mathcal{E}^2(u, t) = \frac{\int_\Omega \lvert \mathcal{N}(u; W)(\cdot, t) - G(u)(\cdot, t)\rvert^2 \, \dx}{\int_\Omega \lvert G(u)(\cdot, t) \rvert^2 \, \dx},
   	\end{equation}
        where as before the latter term is the relative error at a specific time point $t$.

	\subsection{Poisson equation}\label{sec:examples:poisson}
	In this section, we consider the Poisson equation
	\begin{equation}
		\begin{aligned}
			\nabla \cdot \kappa \nabla u &= f \quad \text{ in } \Omega, \\
			u &= g \quad \text{ on } \partial \Omega
		\end{aligned}
		\label{eq:examples:poisson:equation}
	\end{equation}
	with $\Omega = (0, 1)^2$, source term $f : \Omega \rightarrow \mathbb{R}$, solution $u : \Omega \rightarrow \mathbb{R}$, and boundary condition $g : \partial \Omega \rightarrow \mathbb{R}$.
	The spatially varying coefficient $\kappa : \Omega \rightarrow \mathbb{R}$ is assumed unknown.

        Using the hybrid FEM-NN approach we shall investigate (i) how the error in the predicted coefficient $\kappa$ decreases with discretisation error,
        (ii) sensitivity of the neural network to noise in the observations of $u$ and
        (iii) the effect of partial observations of $u$ on the predictions of the coefficient.
	For all experiments, $\kappa(x, y)$ was approximated by a NN $\mathcal{N}(\cdot; W):\Omega\rightarrow\mathbb{R}$ with a single hidden layer
        consisting of 30 neurons and sigmoid activation functions. Note that $\mathcal{N}$ is here  independent of the state.

	\subsubsection{Discretisation error}
	First, we investigate how the prediction error of the neural network depends on the discretisation error of the finite element method.
	We let the analytical solution $u$ and diffusion coefficient $\kappa$ be
	\begin{align*}
	    u(x, y) = \sin(\pi x) \sin(\pi y), \quad \kappa(x, y) = \frac{1}{1 + x^2 + y^2}.
	\end{align*}

        Using gradually refined meshes $\Omega_h$ we train five different neural networks $\mathcal{N}_h \approx \kappa$.
        The architecture of the neural network remained the same for all five meshes.
        On each mesh the spaces $U_h=V_h$ are constructed using continuous piecewise-linear ($P_1$) functions.
	Starting with the same initial weights, the five neural networks are trained for up to 50\,000 L-BFGS iterations with the loss function:
	\begin{align*}
	    L = \int_{\Omega_h} (\hat u_h - u)^2 \, \dx 
	\end{align*}
	where $\hat u_h$ is the predicted solution of \cref{eq:examples:poisson:equation} on $\Omega_h$.
	We note that as the mesh is refined, the integration of $L$ evaluates the true solution at more points and thus more observations are used.
	
\begin{figure}[!htb]
		\centerline{\includegraphics[width=0.5\linewidth]{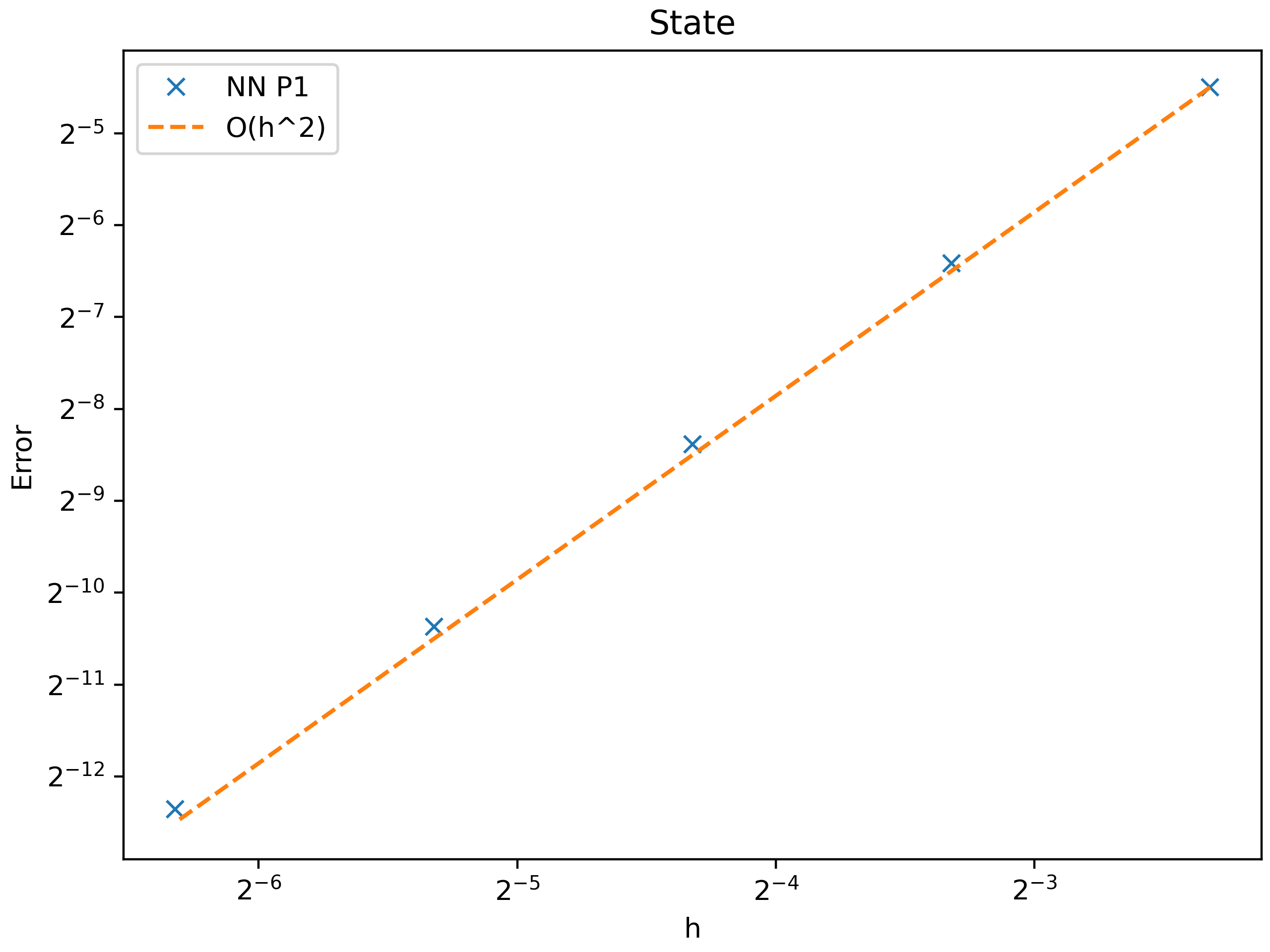}
			\includegraphics[width=0.5\linewidth]{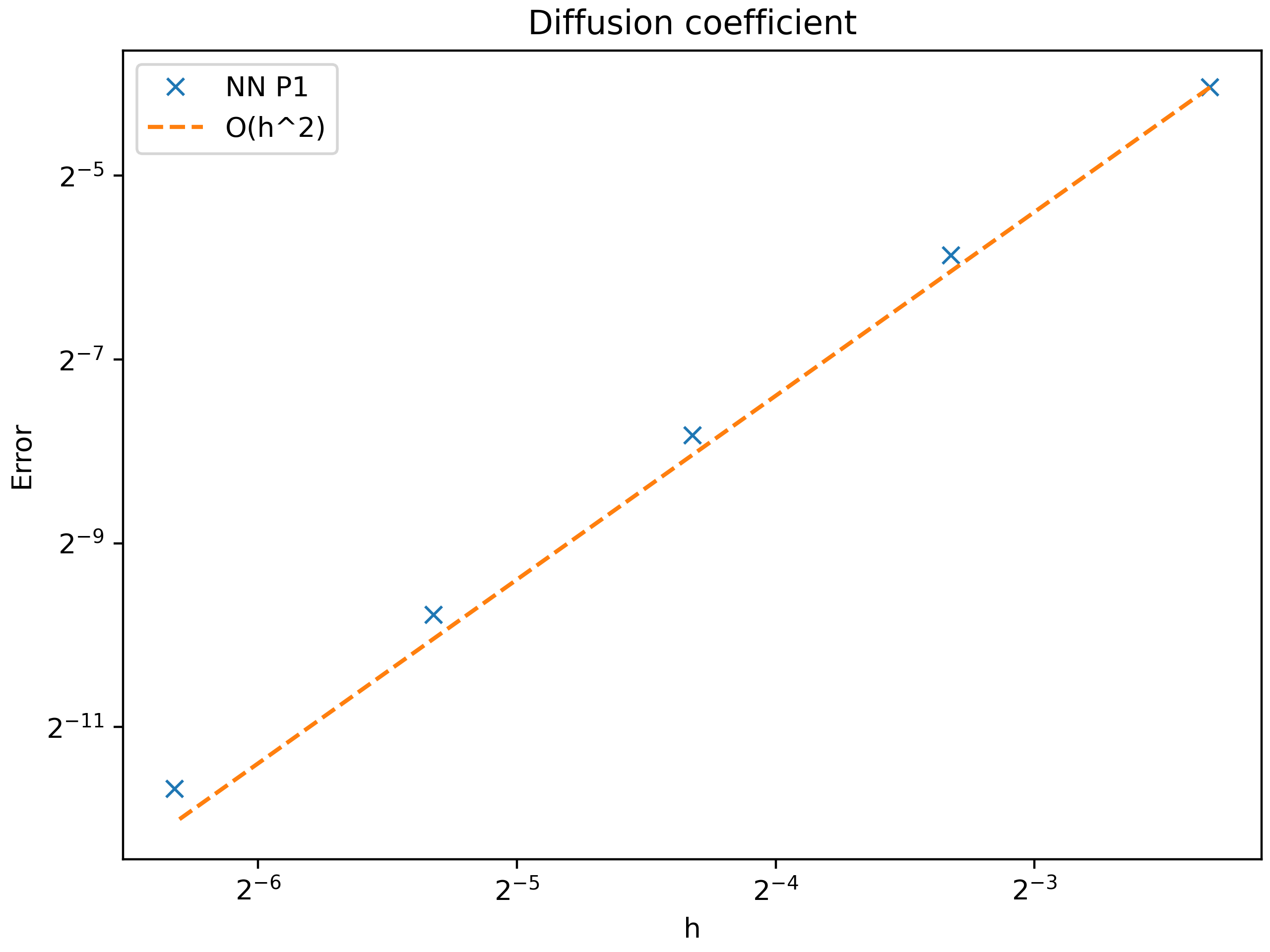}}
		\caption{The prediction errors for different mesh resolutions using $P_1$ elements for the state.
                  Left: the $L_2(\Omega)$-error in the predicted state when trained on a mesh with distance $h$ between vertices. The predicted states are produced by using FEM on a high resolution mesh.
                  Right: the $L_2(\Omega)$-error in the neural network compared to the analytical coefficient $\kappa$.
                }
		\label{fig:examples:poisson:disc:state}
	\end{figure}

	We are interested in how the errors of the trained neural networks influence the predicted solutions. Therefore, in order
        to avoid discretisation errors, we re-evaluated the five hybrid models on a high resolution mesh and piecewise
        second-order polynomials. The error in these high-resolution solutions can be seen in \cref{fig:examples:poisson:disc:state} together with the error in the approximated diffusion coefficient.
	We observe that the error in both the predicted solution and diffusion coefficient decrease quadratically for $P_1$ elements.
        This observation is in agreement with the analysis 
	in \cite{huang_learning_2020} where an upper bound on the approximation error of a NN for a 1-dimensional Poisson
        problem and $P_1$ finite elements
	\begin{equation*}
		\norm{\kappa - \mathcal{N}_h} \leq C_1 \epsilon + C_2 h^2.
	\end{equation*}
	Here $\epsilon$ is the optimisation error such that the objective functional is bounded by $O(\epsilon)$ and $C_1, C_2$ depend
        on the data.
	
	\subsubsection{Noisy observations}\label{sec:noisy_observations}
	In this example we investigate the performance of the FEM-NN trained neural network on noisy data.
	We again consider \cref{eq:examples:poisson:equation} with the analytical solution and diffusion coefficient
	\begin{align*}
	    u(x, y) = \sin(\pi x) \sin(\pi y), \quad \kappa(x, y) = \frac{1}{1 + x^2 + y^2 + (x - 1)^2 + (y - 1)^2}.
	\end{align*}
    The objective functional is defined as
	\begin{align*}
		L = \norm{\hat u - (u + \epsilon)}_{L_2(\Omega)}^2 + \alpha \norm{\mathcal{N}}_{H_0^1(\Omega)}^2,
	\end{align*}
	where $\hat u$ is the predicted solution, $u$ is the analytical solution and $\epsilon$ is a piecewise linear interpolation of pointwise noise.
	That is, at each interior degree of freedom $x_i \in \Omega_h$, we add normally distributed noise $\epsilon(x_i) \sim \mathcal{N}(0, \sigma^2)$ with standard deviation $\sigma = \norm{u}/10 = 0.05$. This results in a signal to noise ratio $\norm{u}/\sigma = 10$ in the $L_2(\Omega)$ norm.
	The last term in the objective functional is a $H_0^1$ regularization term for the neural network with regularization parameter $\alpha \geq 0$.
    
\begin{figure}[!htb]
		\centerline{\includegraphics[width=0.98\linewidth]{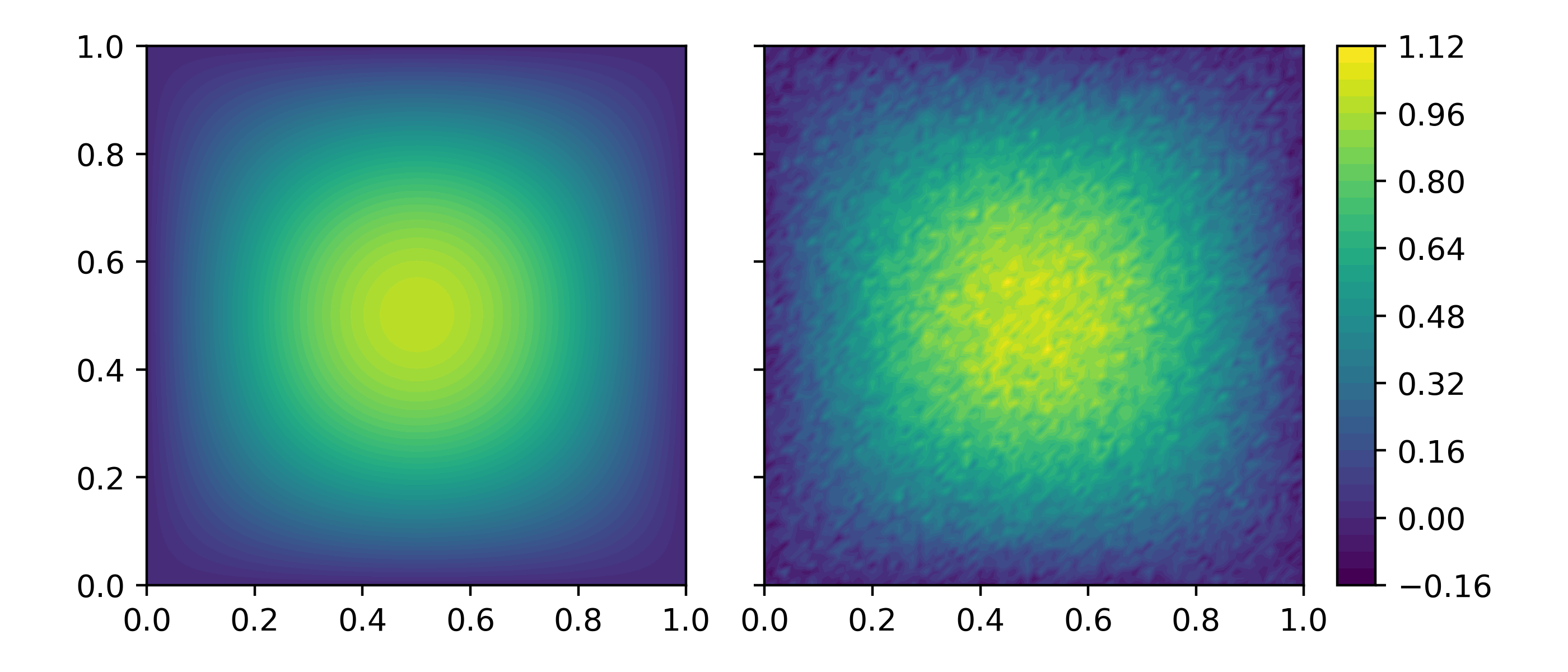}}
		\centerline{\includegraphics[width=0.98\linewidth]{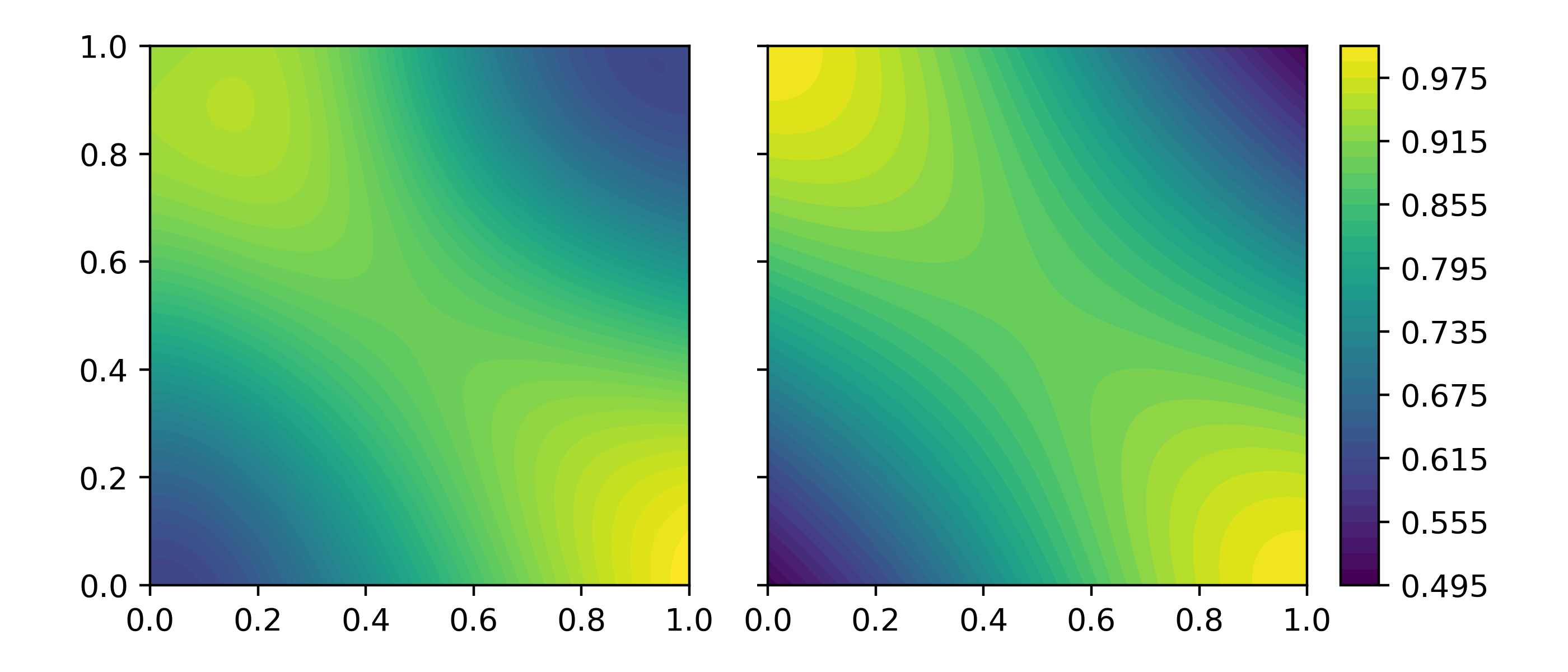}}
		\caption{Poisson problem with noisy data, section \ref{sec:noisy_observations}. Top left: predicated state $\hat u$ after training. The predicted state is visually indistinguishable from the analytical solution $u$.  Top right: the noisy observation $u + \epsilon$ used to train the FEM-NN problem.
		Bottom left: the diffusion coefficient predicted by the neural network. Bottom right: The analytical diffusion coefficient $\kappa$.}
		\label{fig:examples:poisson:noise:net}
	\end{figure}

	The hybrid model was trained for 1\,000  L-BFGS iterations with regularization parameter $\alpha = 10^{-4}$.
	The results are shown in \cref{fig:examples:poisson:noise:net}.
	The FEM-NN model is able to re-construct the unknown diffusion coefficient. The state error was $e_{\Omega} = 0.560\%$, while the error in the diffusion coefficient was $\mathcal{E} = 1.98\%$.
	
\begin{figure}[!htb]
		\centerline{\includegraphics[width=0.78\linewidth]{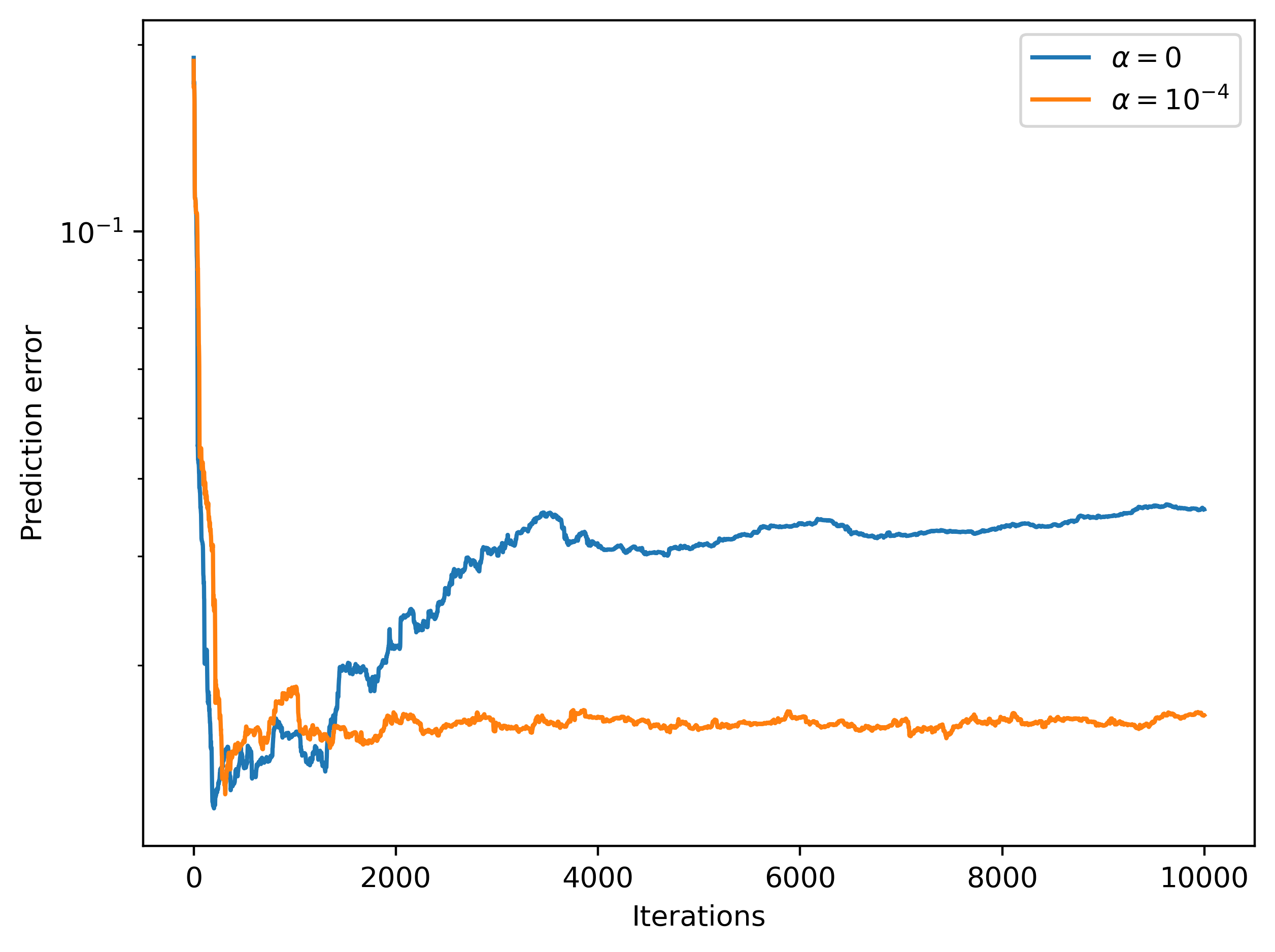}}
		\caption{The prediction error (i.e. $\norm{\kappa - \mathcal{N}}_{L^2(\Omega)}$) of the neural network during training with no regularization (blue line) and with $H_0^1$ regularization (orange line).}
		\label{fig:examples:poisson:noise:convergence}
	\end{figure}

	Berg and Nyström \cite{berg_neural_2017} considered a similar Poisson inverse problem with a different diffusion coefficient and less noise in the observations.
	They compare the performance of the neural network to that of a FEM function represented in the solution space.
	They suggest that the neural network provides a form of implicit regularization that removes the need for $H^1(\Omega)$ regularisation that is required when the unknown diffusion coefficient is represented  as a FEM function.
	We therefore repeated the above experiment with $\alpha = 0$ and $\alpha = 10^{-4}$ and observed the error of the neural network during the training process with 10\,000 L-BFGS iterations. The results are shown in figure \ref{fig:examples:poisson:noise:convergence}.
	Without regularisation, the prediction error drops rapidly during the first 300 iterations. After ca 1\,500 iterations, the prediction errors increase again, caused by overfitting of the neural network to the noise.
	With regularisation, the error drops rapidly during the first 400 iterations, and remains low during the remaining training process.
	Hence, the $H_0^1$ regularization increases the robustness of the training process, increasing resilience to overfitting such that the accuracy of the trained network is not sensitive to the stopping criterion.
	Berg and Nyström \cite{berg_neural_2017} employ a rather sensitive stopping criterion for the optimization algorithm that we observe to stop the optimization loop before the neural network starts to overfit.
	Since the ground truth is realistically unavailable, tuning the stopping criterion is not a reliable strategy for avoiding overfitting.
	The results suggest that also hybrid models require some explicit regularization.
	
	\subsubsection{Partial observations}
	We next consider two variants of \cref{eq:examples:poisson:equation} in which the data are available
        only in parts of the domain's closure $\overline\Omega$. In particular,
        the observations shall be given on a subdomain of $\Omega$ or only on part of the
        boundary. In both cases $P_1$ finite element spaces 
        are used.

        Let us first consider the subdomain $D = \{(x, y) \in \Omega \, : \, |x - 0.5| \leq 0.25, \, |y - 0.5| \leq 0.15 \}$, being a
        rectangle centered at the middle of the unit square with sides of length 0.5 and 0.3,
        inside which we have data on every degree of freedom/vertex of the mesh.
	We let the solution of \cref{eq:examples:poisson:equation} be $u(x, y) = \sin(\pi x) \sin(\pi y)$ and the unknown diffusion coefficient $\kappa(x, y) = 1 + \frac{1}{2} \sin(\pi x) \sin(\pi y)$.
	
\begin{figure}[!htb]
		\centerline{\includegraphics[width=0.5\linewidth]{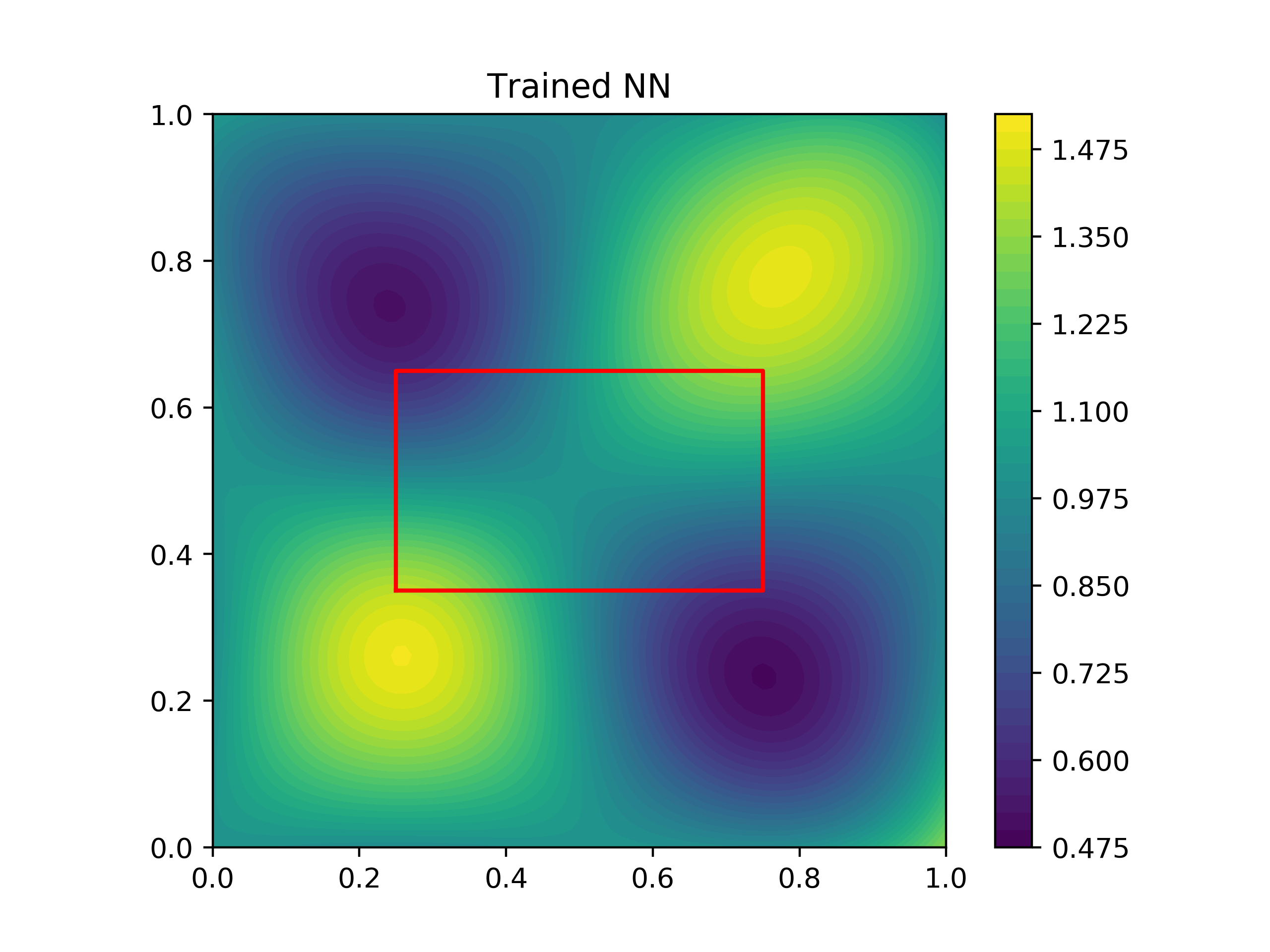}
			\includegraphics[width=0.5\linewidth]{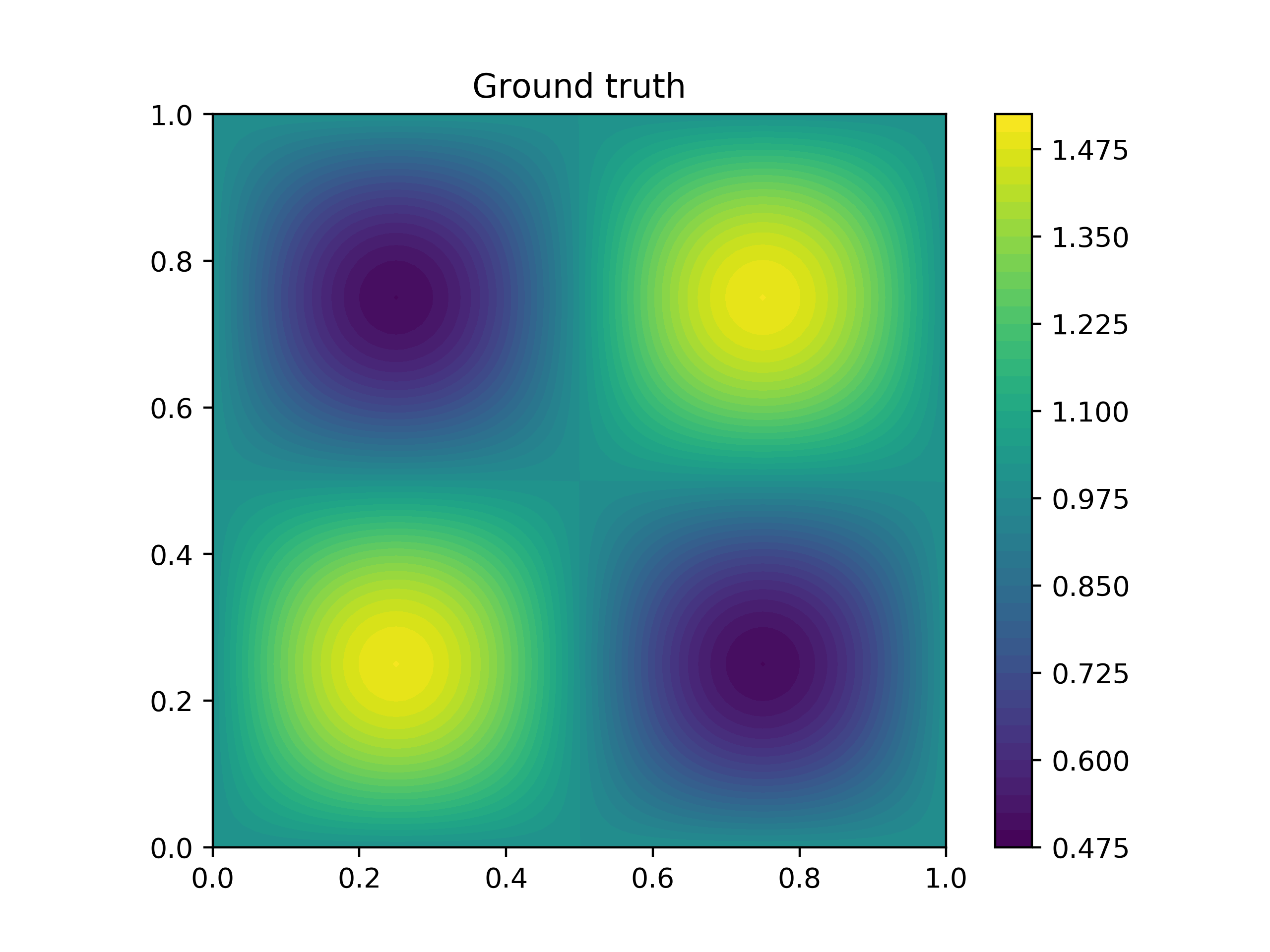}}
		\caption{Left: the trained neural network approximating the diffusion coefficient $\kappa$. The network was only trained on observations from inside the red box. Right: the analytical diffusion coefficient $\kappa$.}
		\label{fig:examples:poisson:partial:nn}
	\end{figure}

	The NN was trained for 10\,000 iterations of L-BFGS and can be seen in \cref{fig:examples:poisson:partial:nn}.
	Although the neural network is only provided observations inside the rectangle,
	it is able to recover an accurate approximation of the analytical diffusion coefficient $\kappa$ also outside the box.
	
	For the second partial observations example we consider the Calderon problem, see \cite{lohner_revisiting_2019}
        and references therein.
	In this problem, there is no source term, $f = 0$, and observations on $u$, $\kappa \frac{\partial u}{\partial n}$
        are given on the boundary $\partial \Omega$. Here, $n$ is the outward directed boundary normal.
        We remark that when $u$ is the voltage potential, $\kappa \frac{\partial u}{\partial n}$ gives
        the current flowing through the boundary and the task of recovering $\kappa$ from
        the voltage/current measurements is of importance in electrical impedance tomography.

	Adopting the example from \cite{lohner_revisiting_2019}, we define the boundary function as:
	\begin{equation*}
		g(x, y) = e^{-\frac{(x-1)^2 + (y - 0.5)^2}{0.5^2}} - e^{-\frac{x^2 + (y - 0.5)^2}{0.5^2}} \quad (x, y) \in \partial \Omega,
	\end{equation*}
	which provides one radial sink term in the middle of the left boundary and a radial source term in the middle of the right boundary.
	The assumed unknown diffusion coefficient is defined as 
	a radial coefficient
	\begin{equation*}
		\kappa(x, y) = 1 + 4 e^{-\frac{(x - 0.5)^2 + (y - 0.5)^2}{0.2^2}}
	\end{equation*}
	that attains the value 5 in the center of the unit square and decreases exponentially closer to the boundary.
	Synthetic data was generated by solving the equation using the FEM with second order polynomials on triangles using the diffusion coefficient $\kappa$. The domain was divided into $2\times 200^2$ triangles. 
	
	In the hybrid FEM-NN approach, the boundary observations on $u$ are used to enforce the Dirichlet boundary condition in the approximated equation, and the derivative information is used in the loss function
	\begin{equation*}
		L = \norm{\mathcal{N} \frac{\partial \hat u}{\partial n} - \kappa \frac{\partial u}{\partial n}}_{L_2(\partial \Omega)}.
	\end{equation*}
	The neural network was first trained for 1\,000 iterations of L-BFGS. At this 
	point, due to slow convergence, we switched 
	to the truncated Newton method (TNC) and performed additional 1\,000 iterations.

\begin{figure}[!htb]
		\centerline{\includegraphics[width=0.5\linewidth]{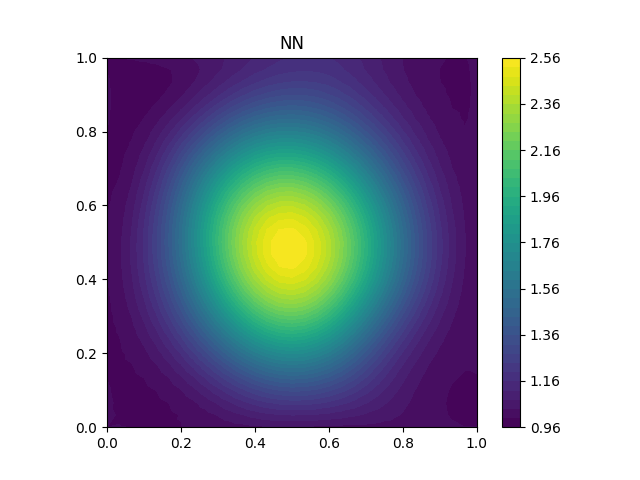}
			\includegraphics[width=0.5\linewidth]{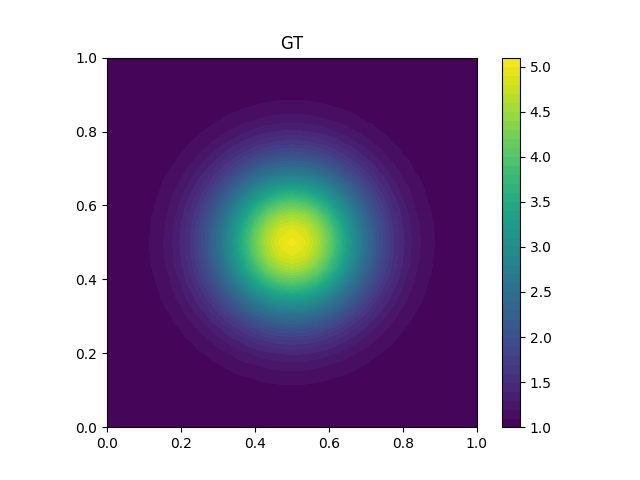}}
		\caption{Left: the trained neural network approximation of the diffusion coefficient $\kappa$. Right: the analytical diffusion coefficient $\kappa$. The neural network has only been trained on observations on the boundary.}
		\label{fig:examples:poisson:calderon:nn}
	\end{figure}

	The resulting neural network approximation is shown in \cref{fig:examples:poisson:calderon:nn}.
	Although the neural network is not able to accurately capture the magnitude of the analytical diffusion coefficient in the center, it is able to capture 
	the radial shape; a significant improvement over the results presented 
	in \cite{lohner_revisiting_2019}. 
	Note that the contribution from the second term in $\kappa$ is on the order of $10^{-3}$ to $10^{-5}$ along the boundary.
	
	\subsection{Comparison to other approaches}
	\label{sec:examples:comparison}
	In this section, we compare the proposed hybrid FEM-NN methodology to physics informed neural networks
        (PINNs) \cite{raissi_physics-informed_2019} and pointwise estimation where the unknown coefficient is represented as a FEM function rather than a NN.
	To this end let us consider an inverse problem for the 1-dimensional heat equation 
	\begin{equation}
		u_t - \nabla \cdot (\kappa \nabla u) = 0 \quad \text{ in } \Omega \times (0, T)
		\label{eq:examples:heat-equation:true}
	\end{equation}
	with $\Omega = (0, 1)$, $T = 0.1$. The diffusion coefficient $\kappa : \Omega \rightarrow \mathbb{R}$ is assumed to be unknown.
	The problem is closed with an initial condition and homogeneous Dirichlet boundary conditions
	\begin{align*}
		u &= 0 \quad \text{ on } \partial \Omega \times (0, T), \\
		u &= u_0 \quad \text{ on } \Omega \times \{0\}
	\end{align*}
	with $u_0(x) = x(x - 1)$.
	Further, we assume $\kappa > 0$ such that the problem remains well-posed.
	
	We apply PINNs and FEM-NN to reconstruct two different choices for the diffusion coefficient
        $\kappa$: a linear case and a discontinuous piecewise constant case.
	Here, the linear case admits a smooth solution $u$ and serves as a basic benchmark for our comparison.
        The discontinuous diffusion coefficient results in a more challenging problem, which, however
        is highly relevant e.g. in biomechanics \cite{valnes2020apparent}. More precisely, the solution
        is not smooth, with jumps in the gradients at the points of the (coefficient) discontinuity.
	
	For each case, synthetic data is generated by solving the equation using a high resolution FEM discretisation in space and explicit first stage implicit 4th order Runge Kutta (ESDIRK43a) as described in \cite{kvaerno_singly_2004}.
	We employed the ESDIRK43a implementation provided in the gryphon-project \cite{skare2012gryphon}.
	The spatial interval is divided into $N_x = 4096$ sub-intervals and the time step was kept fixed at $\Delta t = 10^{-4}$.
        Finally, the set of training observation times is $\mathcal T = \{i \tau \}_{i=1}^{5}$
        where $\tau = 0.02$ is the temporal spacing between (equispaced)
        observation snapshots.

        In the hybrid FEM-NN approach, we replace the unknown term $\kappa$ by a neural network $\mathcal{N}(\cdot; W) : \Omega \rightarrow \mathbb{R}$ which does not depend on $u$ or time. Here the network contains a single hidden layer, hyperbolic tangent activation functions, and weights $W$.
	The hidden layer has 30 neurons. The problem equation then reads
	\begin{equation}
		\begin{aligned}
		\hat u_t - \nabla \cdot (\mathcal{N}(\cdot; W) \nabla \hat u) &= 0 \quad \text{ in } \Omega \times (0, T), \\
		\hat u &= 0 \quad \text{ on } \partial \Omega \times (0, T), \\
		\hat u &= u_0 \quad \text{ on } \Omega \times {0}.
		\end{aligned}
		\label{eq:examples:heat-equation-pinn:approx}
	\end{equation}
	
	Problem \cref{eq:examples:heat-equation-pinn:approx} is solved using FEM in space and Crank-Niholson in time.
	The domain $\Omega$ is discretised into 240 sub-intervals and the $P_1$ finite elements are used. We remark
        that the test functions in $V$ satisfy the (homogeneous) Dirichlet boundary conditions by construction.
	Further, we use $\Delta t = 2.5 \times 10^{-3}$ that  coincides with the observation time points.
	The variational problem formulation with function space $V$ reads: Given $\hat u^0 = u_0$
        for $n=1, 2, \cdots$ find $u^n \in V$ such that
	\begin{equation}
	\begin{aligned}
		\int_\Omega \frac{\hat u^{n+1} - \hat u^n}{\Delta t}\cdot v + \frac{\mathcal{N}(\cdot; W)}{2} ( \nabla \hat u^{n+1} + \nabla \hat u^{n}) \cdot \nabla v \, \dx &= 0 \quad \forall v \in V \\
	\end{aligned}
	\label{eq:examples:heat-equation-pinn:approx:varform}
	\end{equation}
	with $\hat u^n=\hat u(x, n\Delta t)$.
	Hence, the inverse problem formulation in the hybrid FEM-NN model becomes 
	\begin{equation}
	\begin{aligned}
		&\min_{W, \hat u} \sum_{t_i \in \mathcal{T}} \int_\Omega (\hat u(x, t_i) - u(x, t_i))^2 \, \dx
		+ \alpha \int_\Omega (\mathcal{N}(x; W))^2 \, \dx
		\text{ subject to \cref{eq:examples:heat-equation-pinn:approx:varform}}.
		\end{aligned}
		\label{eq:examples:heat-equation-pinn:hybrid-min-problem}
	\end{equation}
	Here $\alpha > 0$ is an $L^2(\Omega)$ regularisation parameter for the neural network.
	This regularisation is added due to the inverse problem being ill-posed, see \cite{shen1999numerical}.
	For our problems we choose a regularisation parameter $\alpha = 10^{-8}$.
        
	We remark that the positivity constraint $\kappa \approx \mathcal{N} > 0$ is only ensured during initialization of the weights $W$ before training.
	Although $\mathcal{N}$ could violate the positivity constraint during training,
	we did not encounter this problem for the results presented below.
	In those cases, one could consider using ReLU or a smooth approximation to ReLU on the output to enforce the constraint
        by construction.
	
	The pointwise estimator approach uses the same formulation of the inverse problem as
        \cref{eq:examples:heat-equation-pinn:hybrid-min-problem}, except that the neural network $\mathcal{N}$ is replaced by an unknown $P_1$ function over the (same) mesh of $\Omega$.
	
	The PINN approach does not use a variational formulation of the problem, and instead incorporates \cref{eq:examples:heat-equation-pinn:approx} as penalty to the minimisation problem.
	The solution of \cref{eq:examples:heat-equation-pinn:approx} is approximated using a neural network $\hat u(x, t) = \mathcal{N}_2(x, t; W_u)$ parameterized by $W_u$. This network consists
        of 3 hidden layers with 16 neurons each and hyperbolic tangent activation functions.
    We use a NN, $\mathcal{N}$, with the same architecture as in the FEM-NN approach to approximate $\kappa$.
	The inverse problem in the PINN case reads
	\begin{equation}
		\min_{W, W_u} \frac{w_{\text{data}}}{\lvert\mathcal{T} \times \mathcal{X}\rvert}\sum_{t_i \in \mathcal{T}} \sum_{x \in \mathcal{X}} (\hat u(x, t_i; W_u) - u(x, t_i))^2 + \frac{\alpha}{\lvert\mathcal{X}\rvert} \sum_{x \in \mathcal{X}} (\mathcal{N}(x; W))^2 + L_{\text{pde}}(\hat u(\cdot, \cdot; W_u)).
		\label{eq:examples:heat-equation-pinn:pinn-min-problem}
	\end{equation}
	Here $L_{\text{pde}}$ measures the residual of \cref{eq:examples:heat-equation-pinn:approx} point-wise at collocation points 
	\begin{align*}
		L_{\text{pde}}(\hat u) = w_{\text{int}} F(\hat u) + w_{\text{bc}} L_{\text{BC}}(\hat u) + w_{\text{ic}} L_{\text{IC}}(\hat u),
	\end{align*}
	where $F$ is the $l^2$ residual of the PDE in the interior
	\begin{align*}
		F(\hat u) = \frac{1}{\lvert\mathcal{X}_\text{int}\rvert}\sum_{(x, t) \in \mathcal{X}_{\text{int}}} (\hat u_t(x, t) - \nabla \cdot \mathcal{N}(x; W) \nabla \hat u(x, t))^2
	\end{align*}
        and
        \[
          L_{\text{BC}}(\hat u) = \frac{1}{\lvert\mathcal{X}_\text{bc}\rvert}\sum_{(x, t) \in \mathcal{X}_{\text{bc}}} (\hat u_t(x, t))^2,\quad
          L_{\text{IC}}(\hat u) = \frac{1}{\lvert\mathcal{X}_\text{ic}\rvert}\sum_{x \in \mathcal{X}_{\text{ic}}} (\hat u_t(x, 0) - u_0(x))^2          
        \]
          enforce the boundary and initial conditions.
          Sets $\mathcal{X}_{\text{int}}, \mathcal{X}_{\text{bc}}, \mathcal{X}_{\text{ic}}$ are the sets of collocation points in the interior, boundary and for the initial condition respectively. $\mathcal{X}$ is the set of spatial data points where the misfit with the data along with the regularization term is evaluated.
	Further, $w_i$, $i=\left\{\text{data}, \text{int}, \text{bc}, \text{ic}\right\}$ are weights balancing the different terms in the loss function.
	If not otherwise stated, we let these weights be $1$.
	Finally, we sample the collocation points in a regular grid with $241$ points in $\Omega$, i.e. the vertices of the finite element mesh, and $40$ points in $[0, T]$.
	This coincides with the degrees of freedom for the numerical PDE solution $\hat u$ over $\overline\Omega \times [0, T]$ used in the hybrid FEM-NN method.
	We let the set of spatial data points at each observation, $\mathcal{X}$, be all 241 collocation points.
	This means that the data function $u$ is interpolated at the same points for both the hybrid FEM-NN method and the PINNs method.
	However, it should be noted that the method of incorporating these data points is slightly different in the FEM-NN method, where due to the $L^2$ integral, an additional linear interpolation between the spatial data points occur.

        We note that the form of PINNs considered further is that from the original paper
        \cite{raissi_physics-informed_2019}. However, the method has seen a rapid development with
        a number of PINN variants (e.g. \cite{jagtap2020conservative, jagtap2020extended}),
        training methods \cite{wang2020understanding} or architectures \cite{jagtap2020adaptive} proposed,
        which would likely improve its performance on our benchmark problems. We chose the
        original implementation simply as a point of reference.

	\subsubsection{Linear diffusion coefficient}
	Let the diffusion coefficient be defined as
	\begin{align*}
		\kappa(x) = 2x + 1.
	\end{align*}
	In the hybrid approach, \cref{eq:examples:heat-equation-pinn:hybrid-min-problem} is optimised using SciPy L-BFGS \cite{10.1145/279232.279236, 2020SciPy-NMeth} for 100 iterations.
	The pointwise estimator is trained for 1\,000 iterations of L-BFGS.
	For the PINN approach, \cref{eq:examples:heat-equation-pinn:pinn-min-problem} is optimised using 100\,000 iterations of L-BFGS.
	The PINN optimisation approach requires significantly more training iterations since two neural networks are trained simultaneously: the unknown diffusion coefficient for the inverse problem and the solution of the PDE given the current estimate of the diffusion coefficient.
	However, the forward and backpropagation of the PINN model is typically faster than for the hybrid model.

	In addition, we trained a NN with the PINNs approach where we adjusted the weights that balance the loss terms.
	Specifically, we manually tuned the weights such that the data term is given greater attention:
	$w_\text{data} = 10^3, w_\text{ic} = 10^2, w_\text{bc} = 10^1$, and $w_\text{int} = 10^0$.
	This network is trained for 10\,000 iterations of L-BFGS.
	
\begin{figure}[!htb]
		\centerline{\includegraphics[width=0.5\linewidth]{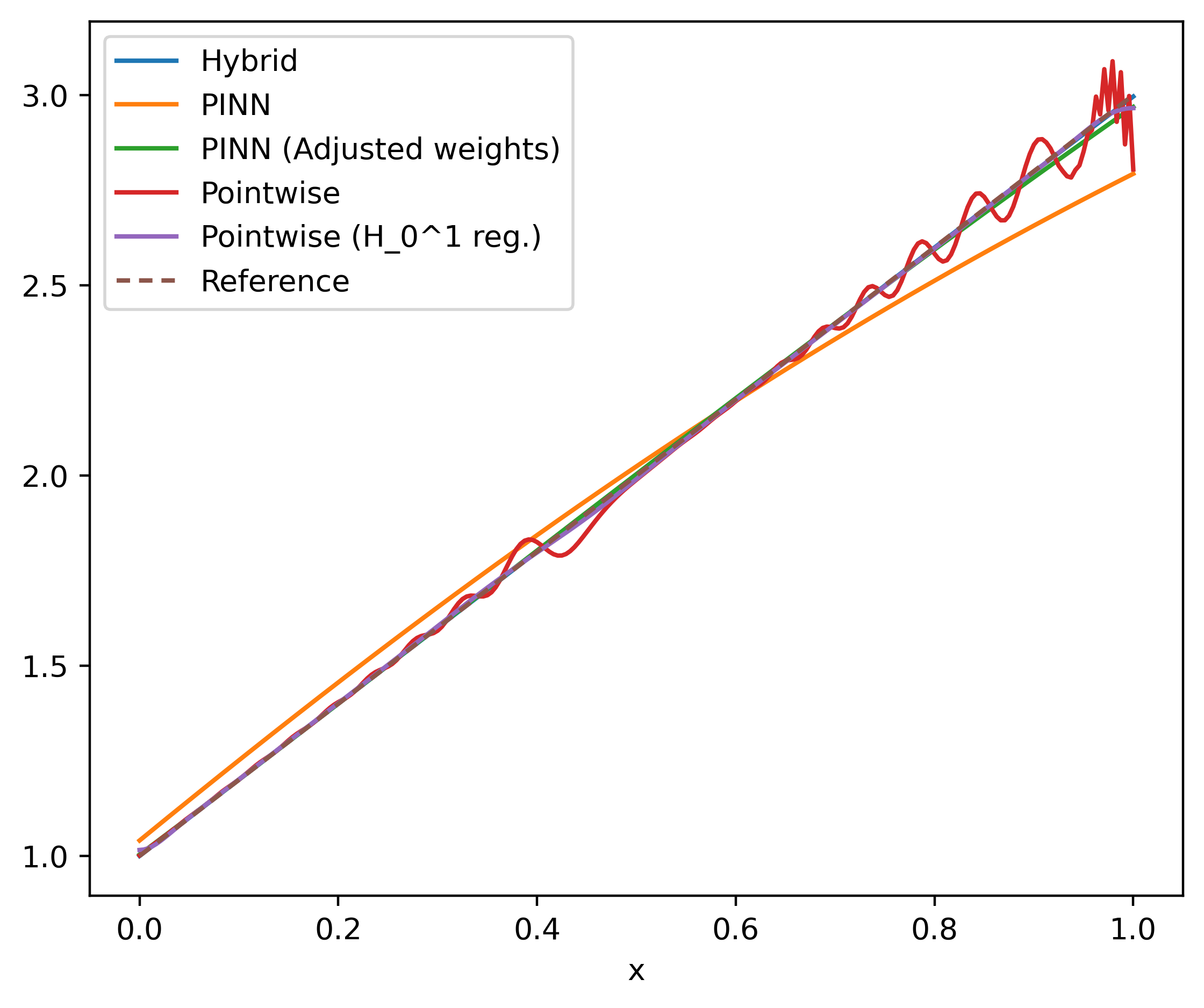}
		\includegraphics[width=0.5\linewidth]{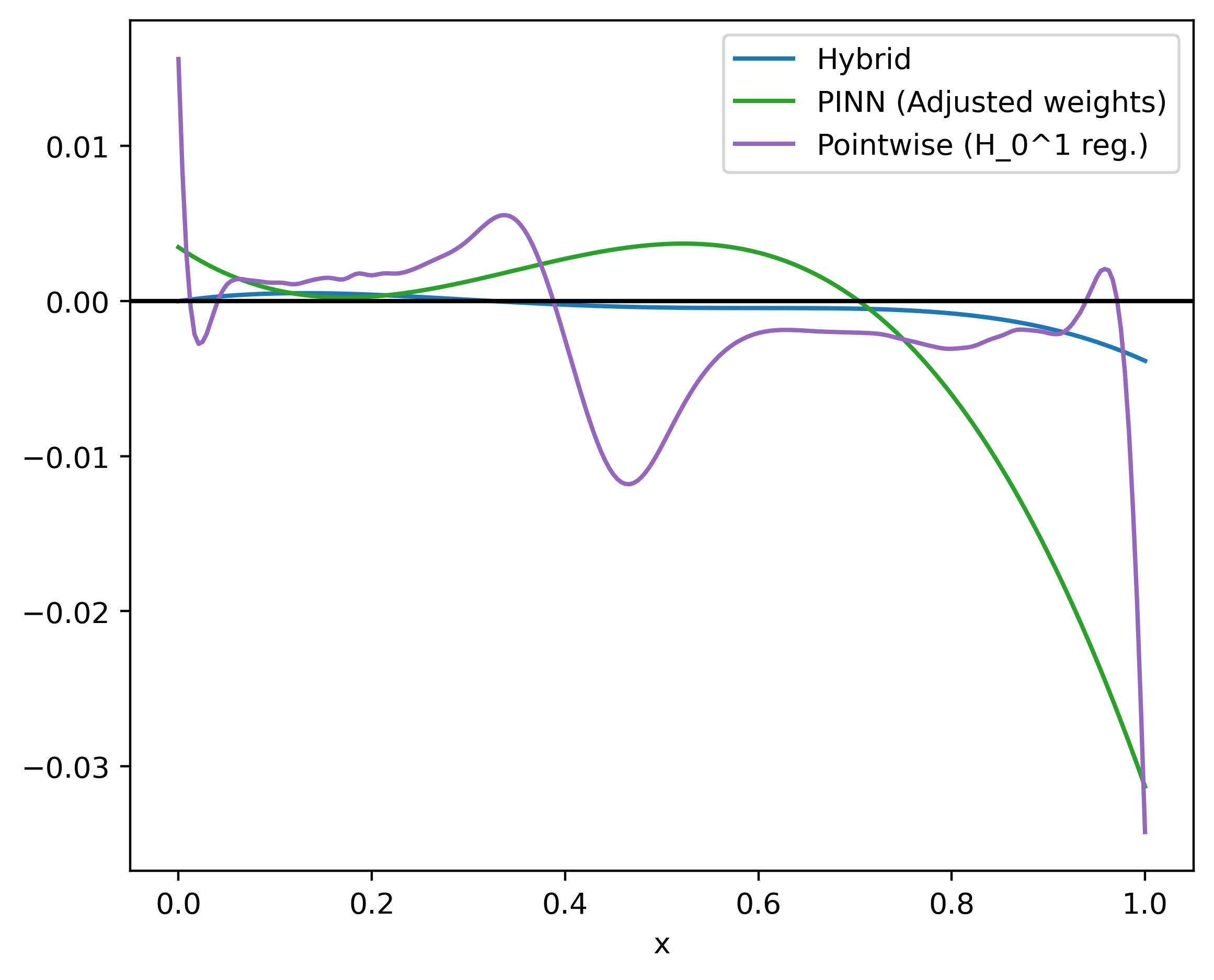}}
		\caption{Left: the trained neural networks as well as the pointwise estimators of the linear diffusion coefficient $\kappa$. The hybrid method, PINNs with adjusted weights, and the pointwise estimator with $H_0^1$ regularisation are able to accurately approximate $\kappa$. Right: the error $\mathcal{N} - \kappa$ for the trained neural networks from the hybrid method and the PINN with adjusted weights, in addition to the pointwise estimator with $H_0^1$ regularisation.}
		\label{fig:examples:heat-equation-pinn:pinn-kappa1}
	\end{figure}

	The trained networks are plotted in \cref{fig:examples:heat-equation-pinn:pinn-kappa1}.
	We observe that both the hybrid and PINN method are able to recover an accurate approximation of the unknown diffusion coefficient $\kappa$.
	However, even after 100\,000 iterations, the PINN trained network with default weights has quite a significant error compared to the true $\kappa$, whereas this error is significantly reduced when using the adjusted weights, even though it is trained with 10 times fewer iterations. The PINNs performance could likely be even further improved by adapting the weights throughout the training process as suggested in \cite{wang2020understanding}. 

	The wall clock run times for the training were a bit under 2 minutes for the hybrid model, roughly 8 minutes for 10\,000 iterations and a bit over 1 hour for 100\,000 iterations of training the PINNs.
    All training is performed on a AMD Ryzen Threadripper 3970X CPU on a single thread, on a machine with 128 GB RAM.
	Thus, the hybrid approach outperforms the PINN method on a CPU architecture, both in terms of training time and in accuracy of the solution.
	This advantage comes mainly from the PINN approach involving two neural networks that balance each other during training.
	In other words, the data terms and PDE term in \cref{eq:examples:heat-equation-pinn:pinn-min-problem} can only decrease with small decrements because a decrease in one term will likely cause an increase in the other term.
    We remark that the use of adaptive activation functions \cite{jagtap2020adaptive} or adaptive
        learning rates \cite{wang2020understanding} in PINNs could speed up convergence.
	
	Further, we observe that in \cref{fig:examples:heat-equation-pinn:pinn-kappa1} the pointwise estimator is not able to accurately recover the diffusion coefficient.
	Although the error in the solution $u$ is comparable to the hybrid method,
	the pointwise estimator produces a non-smooth solution to the inverse problem.
	This is in line with previous results in \cite{berg_neural_2017, fan_solving_2020}, where it is shown that the neural network can act as regularisation and produces smooth approximations.
	It should be noted that the number of degrees of freedom used for the pointwise estimator is significantly greater than the number of weights in the neural network: the pointwise estimator has 241 degrees of freedom while the neural network has 91.
	
	Since $L^2$-regularisation results in an oscillating pointwise estimator, we attempt to impose higher regularity on the estimated solution by replacing the $L^2$-regularisation with $H_0^1$-regularisation:
	\begin{align*}
		\alpha \int_\Omega \left| \nabla \mathcal{N}(x; W) \right|^2 \, \dx,
	\end{align*}
	where $\mathcal{N}$ is the pointwise estimator and $\alpha$ is again chosen to be $\alpha = 10^{-8}$.
	The resulting estimator can be seen in \cref{fig:examples:heat-equation-pinn:pinn-kappa1}.
	It accurately reconstructs the linear diffusion coefficient, and is comparable to the hybrid and PINN method.
	Furthermore, the pointwise estimator converges in less than 100 iterations of L-BFGS, outperforming both PINN and the hybrid method in terms of time spent training.
	However, it still has slightly higher error than the hybrid method.
	This error could likely be further reduced by tuning the regularisation parameter $\alpha$.

	\subsubsection{Discontinuous diffusion coefficient}
	Let the diffusion coefficient be piecewise constant
	\begin{align*}
		\kappa(x) = \begin{cases}
		2 & |x - 0.5| \leq 0.25 \\
		1 & \text{ otherwise.}
		\end{cases}.
	\end{align*}
	In the forward problem classical PINNs cannot capture the discontinuity of $\kappa$ in the PDE residual, with the derivative of $\kappa$ being treated as being everywhere $0$.
	Thus, when solving equation \cref{eq:examples:heat-equation:true} with a classical PINNs formulation, the result is instead the solution to $\hat u_t - \kappa \nabla \cdot \nabla \hat u = 0$.
	One way to mitigate this is to use a weak formulation of the PDE with variational PINNs (VPINNs) \cite{kharazmi_variational_2019}, or by using domain decomposition through extended PINNs (XPINNs) \cite{jagtap2020extended}.
    In XPINNs,  $u$ and $\kappa$ would be represented by several ``subdomain'' networks coupled across the interfaces,
    some of which could be the points of discontinuity, thus making it possible to capture the non-smooth coefficient. The downside is that the number and location of discontinuities must be  known a priori.
	For this reason, we continue with the classical PINNs approach and instead formulate the equivalent mixed problem
	\begin{align*}
		\hat u_t - \nabla \cdot \phi &= 0 &\quad\mbox{ in }\Omega,\\
		\phi - \kappa \nabla \hat u &= 0 &\quad\mbox{ in }\Omega,
	\end{align*}
	where the newly introduced flux variable $\phi$ is approximated by a third neural network $\mathcal{N}_3(x, t; W_\phi)$. This neural network has the same structure as the PDE solution network $\mathcal{N}_2$.
	Note that unlike with the primal formulation
	in the mixed problem all the derivatives required are well 
	defined everywhere inside the domain. Thus standard PINNs approach can be readily applied.
	This formulation results in the following PINN PDE residual
	\begin{align*}
		F(\hat u) = \frac{1}{\lvert \mathcal{X}_{\text{int}} \rvert}\sum_{(x, t) \in \mathcal{X}_{\text{int}}} \left[ (\hat u_t(x, t) - \nabla \cdot \mathcal{N}_3(x, t; W_\phi)^2 + (\mathcal{N}_3(x, t; W_\phi) - \mathcal{N}(x; W) \nabla \hat u)^2 \right].
	\end{align*}
	
\begin{figure}[!htb]
		\centerline{\includegraphics[width=0.4572\linewidth]{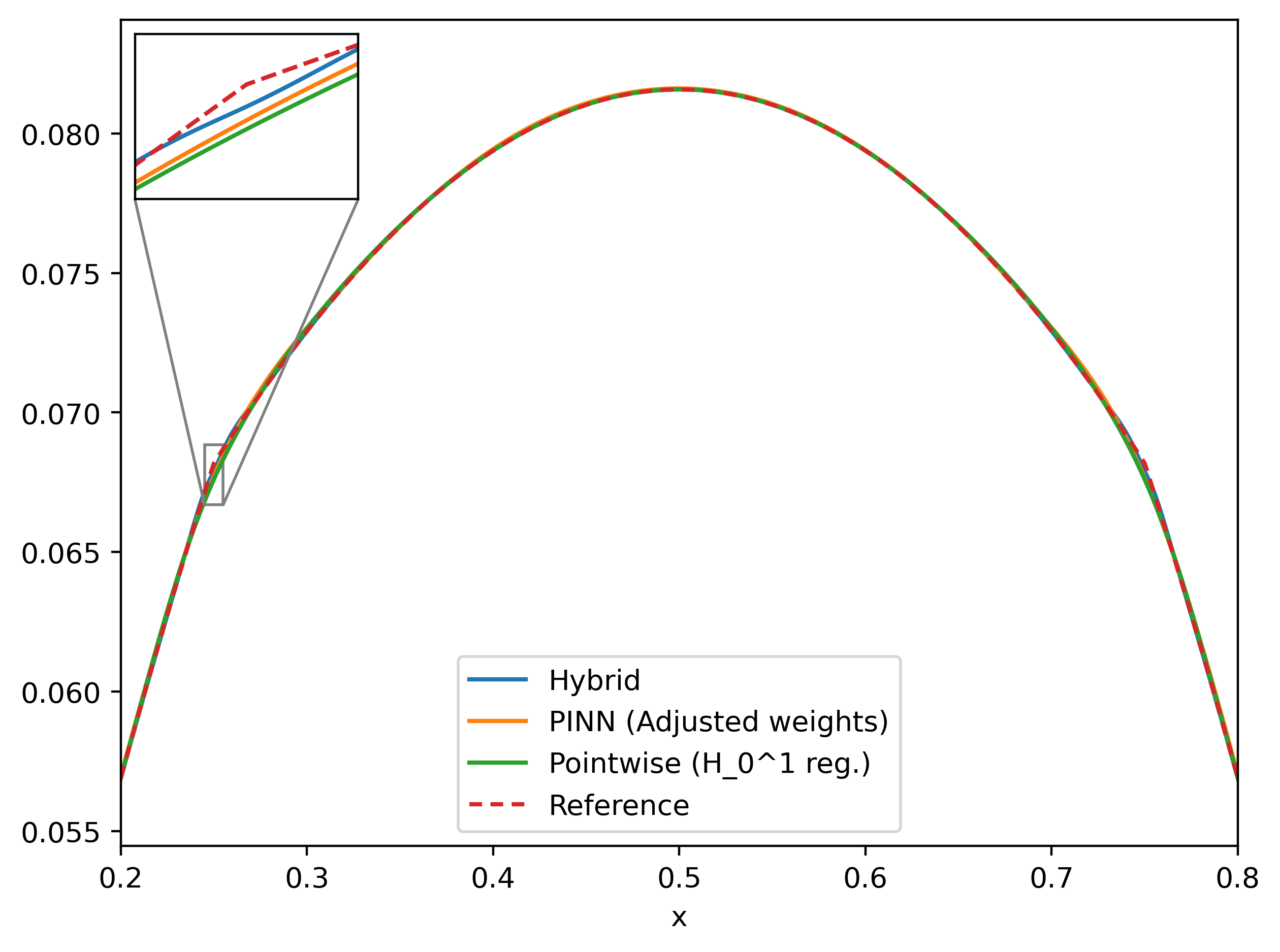}\includegraphics[width=0.4348\linewidth]{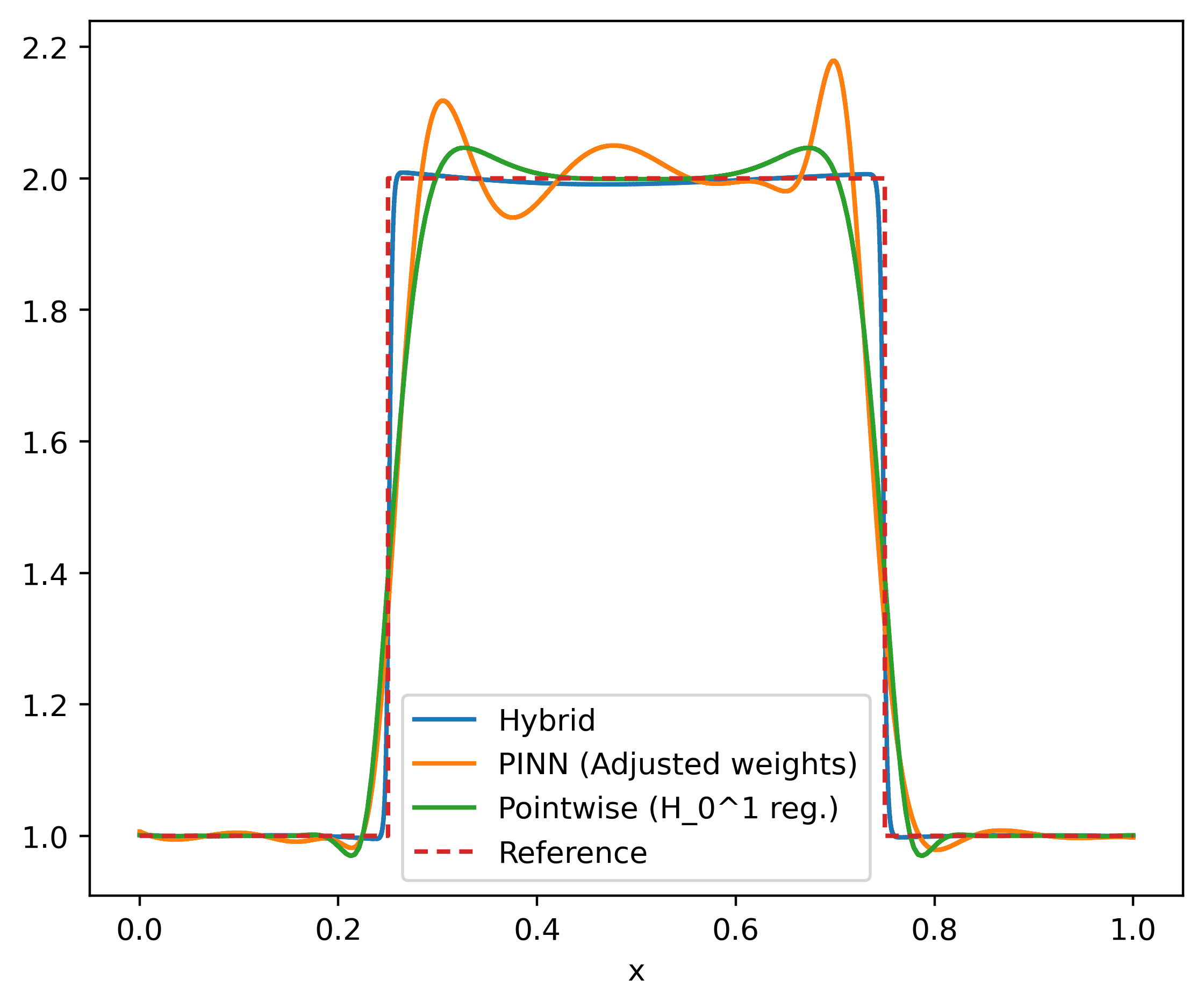}}
		\caption{Left: the approximated solution at $t = 0.1$ of \cref{eq:examples:heat-equation-pinn:approx} with the trained neural networks as well as the pointwise estimator. The plot is limited to the interval $[0.2, 0.8]$ to highlight the differences in the solutions. Right: the approximations of $\kappa$ by the trained neural networks and the pointwise estimator.}
		\label{fig:examples:heat-equation-pinn:pinn-kappa2}
	\end{figure}

	We trained several neural networks with different initial weights for both the hybrid and PINNs method, and selected the best performing network for each.
	The network trained with the hybrid approach was trained for 192 L-BFGS iterations until there was no improvement in the objective functional, followed by 1\,517 TNC iterations (latter required 10\,000 functional evaluations).
	The pointwise estimator was trained for 1\,000 L-BFGS iterations.
	The PINNs network was trained for 100\,000 L-BFGS iterations, followed by 4\,309 TNC iterations (latter required 100\,000 functional evaluations).
	We used the same adjusted weights in the loss function as we used for the linear diffusion coefficient case.
	The trained neural networks, together with the corresponding predicted solutions $\hat u$ at $t = 0.1$ of \cref{eq:examples:heat-equation-pinn:approx}, using the same discretisation scheme as was used to generate the synthetic data, can be seen in \cref{fig:examples:heat-equation-pinn:pinn-kappa2}.
	The predicted solutions are only plotted over the interval $[0.2, 0.8]$ to highlight the differences, since the solutions match the reference solution well outside of this interval.
	
	The hybrid model is able to visually recover the discontinuities in $\kappa$. However, the neural networks with the chosen architecture are smooth functions, and thus cannot represent real discontinuities by design. This can be seen in \cref{fig:examples:heat-equation-pinn:pinn-kappa2} as there is no kink at the discontinuity of $\kappa$ in the predicted solutions $\hat u$ for the hybrid method.
	
	The PINN method has difficulties to recover the unknown diffusion coefficient. The PINNs solution has similarities to what the hybrid model produces in early phases of training, indicating that PINNs could attain better results if trained further. Indeed, we did see some improvements in the PINNs case if we increase the number of TNC evaluations to 1\,000\,000 at the cost of significantly higher computational expense.
	As in the previous experiment there might be benefits to adapt the functional weights $w_i$, $i=\left\{\text{data}, \text{int}, \text{bc}, \text{ic}\right\}$ during training such as suggested in \cite{wang2020understanding}, which was not done here for simplicity.
	
\begin{figure}[!htb]
		\centerline{\includegraphics[width=0.5\linewidth]{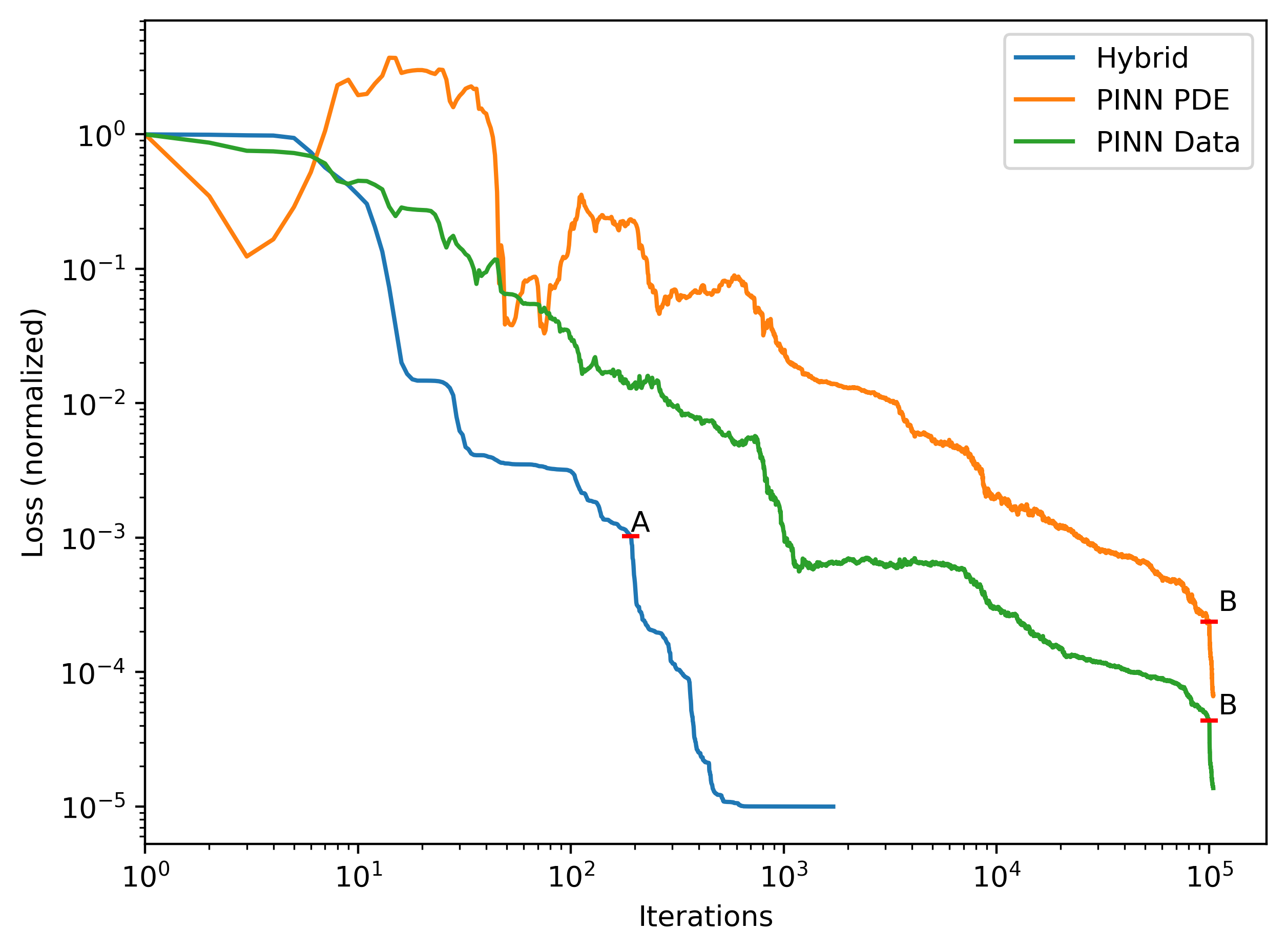}}
		\caption{The normalized loss functions over number of training iterations in a log-log plot. The losses are normalized such that the first iteration is 1. Beyond points A and B, we switch the training of the hybrid method (A) and PINNs method (B) from L-BFGS to the truncated Newton method. It is important to note that few iterations of truncated Newton involve many function evaluations due to inner conjugate gradient iterations.}
		\label{fig:examples:heat-equation-pinn:loss-kappa2}
	\end{figure}

	In \cref{fig:examples:heat-equation-pinn:loss-kappa2} the hybrid and PINN loss functions are plotted along log scaled axes.
	While the PINN data loss, which includes regularisation, attempts to measure the same misfit as the hybrid loss, the two quantities are not directly comparable as one is an integral and the other is the mean over the collocation points.
	Additionally, the reported loss for the PINN is based on the neural network approximation of the state, which might not reflect the actual solution with the current estimate for the diffusion coefficient $\kappa$.
	
	Both the hybrid and PINNs method eventually stagnate with L-BFGS, seeing no improvement beyond the points A and B in \cref{fig:examples:heat-equation-pinn:loss-kappa2}.
	However, switching to TNC allowed the neural networks to improve further, and these TNC iterations were crucial to obtain accurate approximations of $\kappa$.
	
	The pointwise estimator predicts a very irregular estimate of $\kappa$ when using $L^2$-regularisation.
	Thus, as before, we train the pointwise estimator using $H_0^1$-regularisation.
	The resulting estimate of the discontinuous coefficient $\kappa$ can be seen to the right in \cref{fig:examples:heat-equation-pinn:pinn-kappa2}. 
	Again the pointwise estimator with $H_0^1(\Omega)$-regularisation produces a fairly good estimate of the unknown diffusion coefficient $\kappa$.
	While the error is still slightly higher than the hybrid method, the pointwise method converges much faster: at around 2 minutes versus 2 hours for the hybrid method.
	
	\subsection{Advection-diffusion equation}
	\label{sec:examples:advection}
	Although pointwise estimators might be a viable or even better alternative to the hybrid method with an unknown
        spatially varying coefficient like in \cref{eq:examples:heat-equation-pinn:hybrid-min-problem}, pointwise estimation
        is not suited for inverse problems where the unknown term is a differential operator or depends on the state in a
        nonlinear way. We shall now address this type of learning with a hybrid approach.        
        
	We demonstrate the case where the neural network recovers a differential 
	operator. We consider the advection-diffusion equation given by
	\begin{equation}
		u_t - D \Delta u + \nabla \cdot (\omega u) = 0 \quad \text{ in } \Omega \times (0, T),
		\label{eq:examples:advection-diffusion:original}
	\end{equation}
        where $\Omega\subset\mathbb{R}^d$, $D$ is a diffustion constant and
        $\omega:\Omega\rightarrow\mathbb{R}^d$ is the advection velocity. Boundary
        and initial conditions completing \cref{eq:examples:advection-diffusion:original}
        shall be specified later.

	Treating the advection term as unknown we consider a neural 
	network dependent on $u$, the state gradient and the spatial 
	coordinate, i.e. $\mathcal{N}(u; W)(x, t)\equiv\mathcal{N}(x, u(x, t),\nabla u(x, t); W)$ for which
	the learning problem solves the following equation
	\begin{equation}
		\hat u_t - D \Delta \hat u + \mathcal{N}(u; W) = 0 \quad \text{ in } \Omega \times (0, T)
		\label{eq:examples:advection-diffusion:approx}
	\end{equation}
	Note that no information about the advection velocity is available a priori.

        Letting $\Omega=(0, 1)^2$, synthetic data was obtained by solving \cref{eq:examples:advection-diffusion:original} 
        with $D = 0.1$, $\omega(x, y) = [\sin(\pi x) \cos(\pi y), \allowbreak -\cos(\pi x) \sin(\pi y)]$, and $T = 5.0$ with
        initial condition $u(\cdot, 0) \equiv 0$ 
        using $P_1$ elements on a grid with $2\times100\times100$ elements. On the left boundary $x = 0$
        we prescribed Dirichlet boundary conditions $u = \sin^2(t)$ while on the remaining
        boundaries Neumann boundary conditions $\nabla u \cdot n = 0$ were imposed.
        For the temporal discretisation we use the 4th order implicit Runge Kutta with an
        explicit first stage \cite{kvaerno_singly_2004},
        as implemented in the gryphon-project \cite{skare2012gryphon}, with time step $\dt = 0.01$.
	
	For training, the generated synthetic data was sub-sampled in time with a uniform time-step of $0.2$,
        resulting in 26 observations including the initial condition.
	The set of time points for these training samples is denoted by $\mathcal{T}$.
	\Cref{eq:examples:advection-diffusion:approx} was solved using $P_1$ elements on a mesh with $2\times 30 \times 30$
        triangular elements and the Crank-Nicolson scheme with time-step $\dt = 0.1$.
	As the initial and boundary conditions, and diffusion coefficient $D$ are assumed known
        the learning problem becomes
	\begin{equation}
		\min_{W, \hat u} \sum_{t_i \in \mathcal{T}} \norm{\hat u(\cdot, t_i) - u(\cdot, t_i)}_{L^2(\Omega)}^2
		\text{ subject to \cref{eq:examples:advection-diffusion:approx}.}
	\label{eq:examples:advection_diffusion:minimization-problem}
	\end{equation}
        
\begin{figure}[!ht]
		\centering
		\centerline{\includegraphics[width=0.5\linewidth]{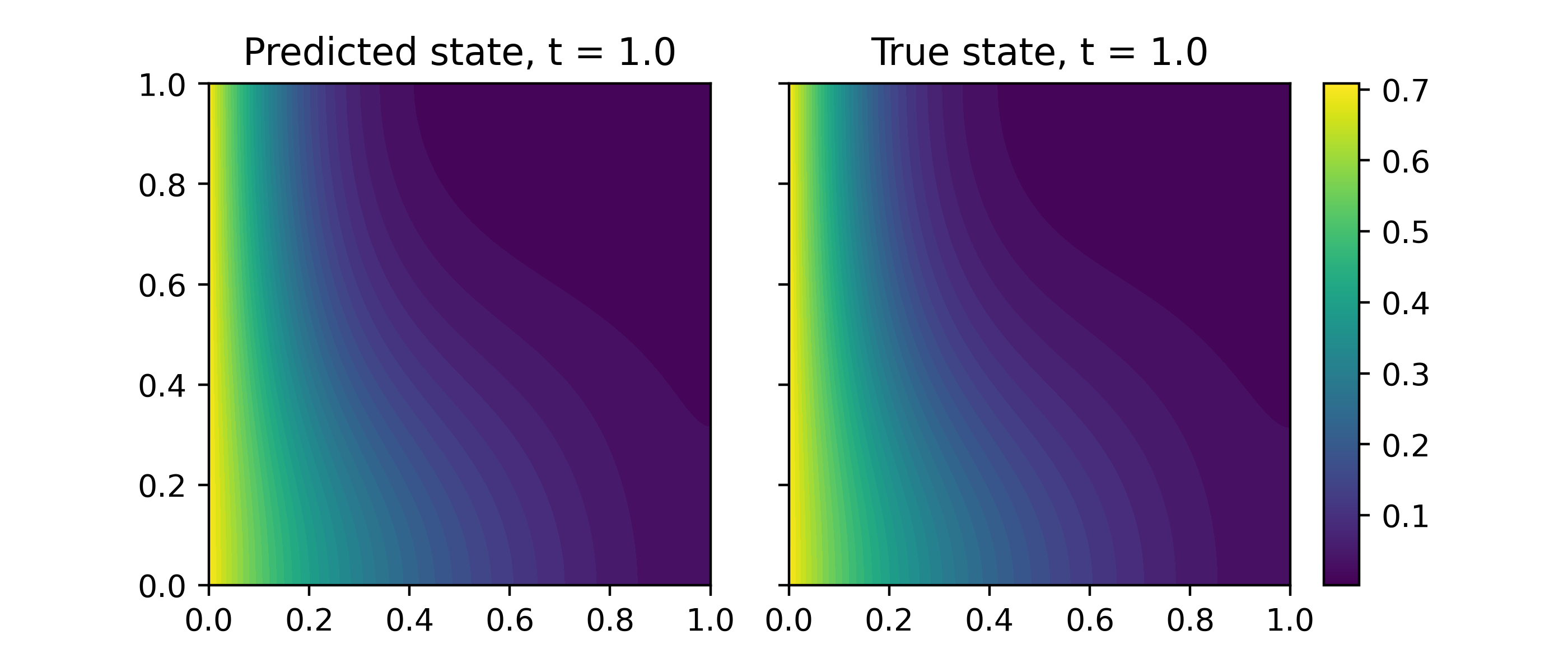}
		\includegraphics[width=0.5\linewidth]{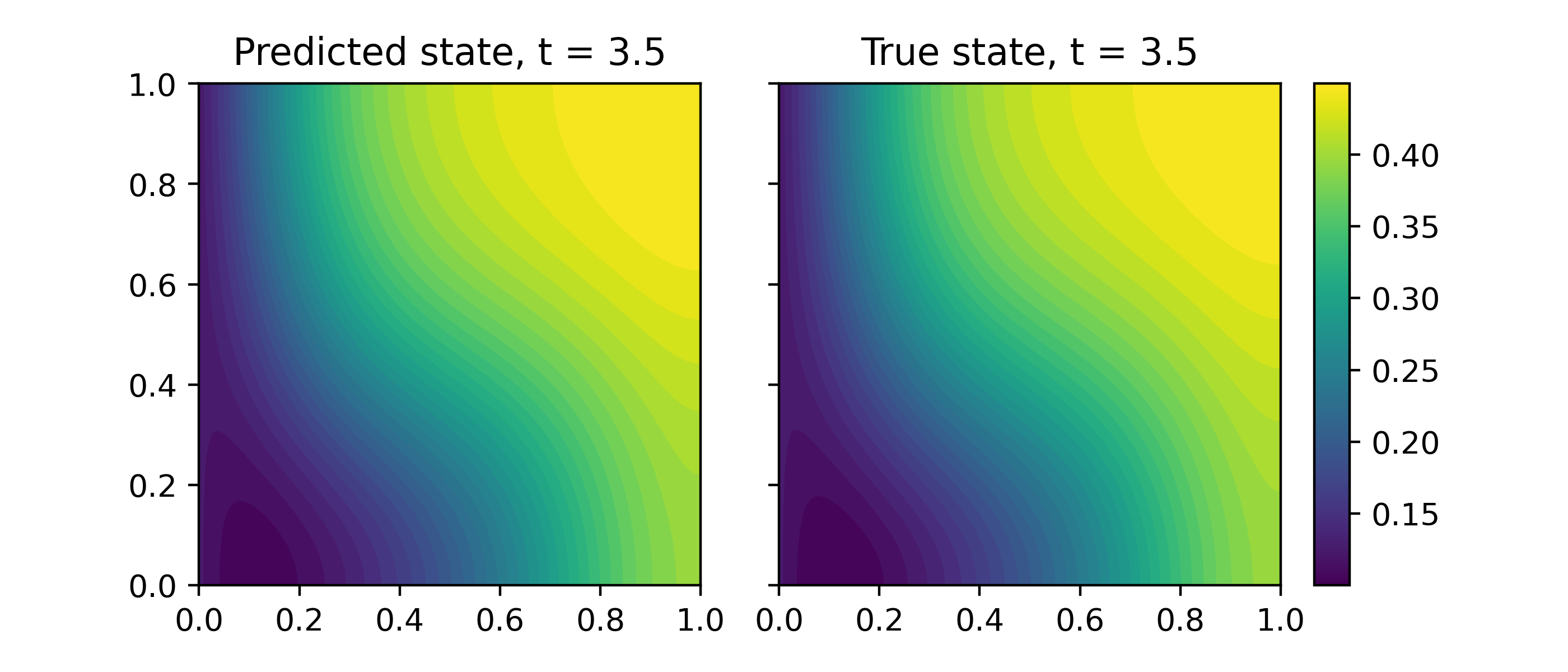}}
		\caption{Left: the predicted solution of \cref{eq:examples:advection-diffusion:approx} at $t = 1.0$ with the fully trained NN,
			and the true solution at $t = 1.0$ based on \cref{eq:examples:advection-diffusion:original}.
			The relative error, $e_{\Omega}(1.0)$, is 0.48\%.
			Right: the predicted solution of \cref{eq:examples:advection-diffusion:approx} at $t = 3.5$ with the fully trained NN,
			and the true solution at $t = 3.5$ based on \cref{eq:examples:advection-diffusion:original}.
			The relative error, $e_{\Omega}(3.5)$, is 0.22\%}
		\label{fig:examples:advection-diffusion:state:1.0+3.5}
	\end{figure}

	We let $\mathcal{N}$ be a neural network with one hidden layer, consisting of 30 neurons, with $\tanh$ activation functions. The network is trained using L-BFGS and TNC.
	The trained neural network is then used to solve \cref{eq:examples:advection-diffusion:approx} using the same
        discretisation scheme as was used to generate the synthetic data.
	\Cref{fig:examples:advection-diffusion:state:1.0+3.5} shows the resulting predicted solution at $t = 1.0$ and $t = 3.5$.
	When evaluating the trained neural network with the ground truth discretisation scheme, the computed errors are $e_{\Omega \times \mathcal{T}}(\hat u) = 0.27\%$ in the predicted state and $\mathcal{E}_{\Omega \times \mathcal{T}}(\mathcal{N}) = 13.6\%$ in the predicted advection term.
	For reference, the relative error in the predicted state is $e_{\Omega \times \mathcal{T}} = 30.6\%$ when only solving the PDE part (i.e. with $\mathcal{N} = 0$).

\begin{figure}
		\centering
		\includegraphics[width=0.5\linewidth]{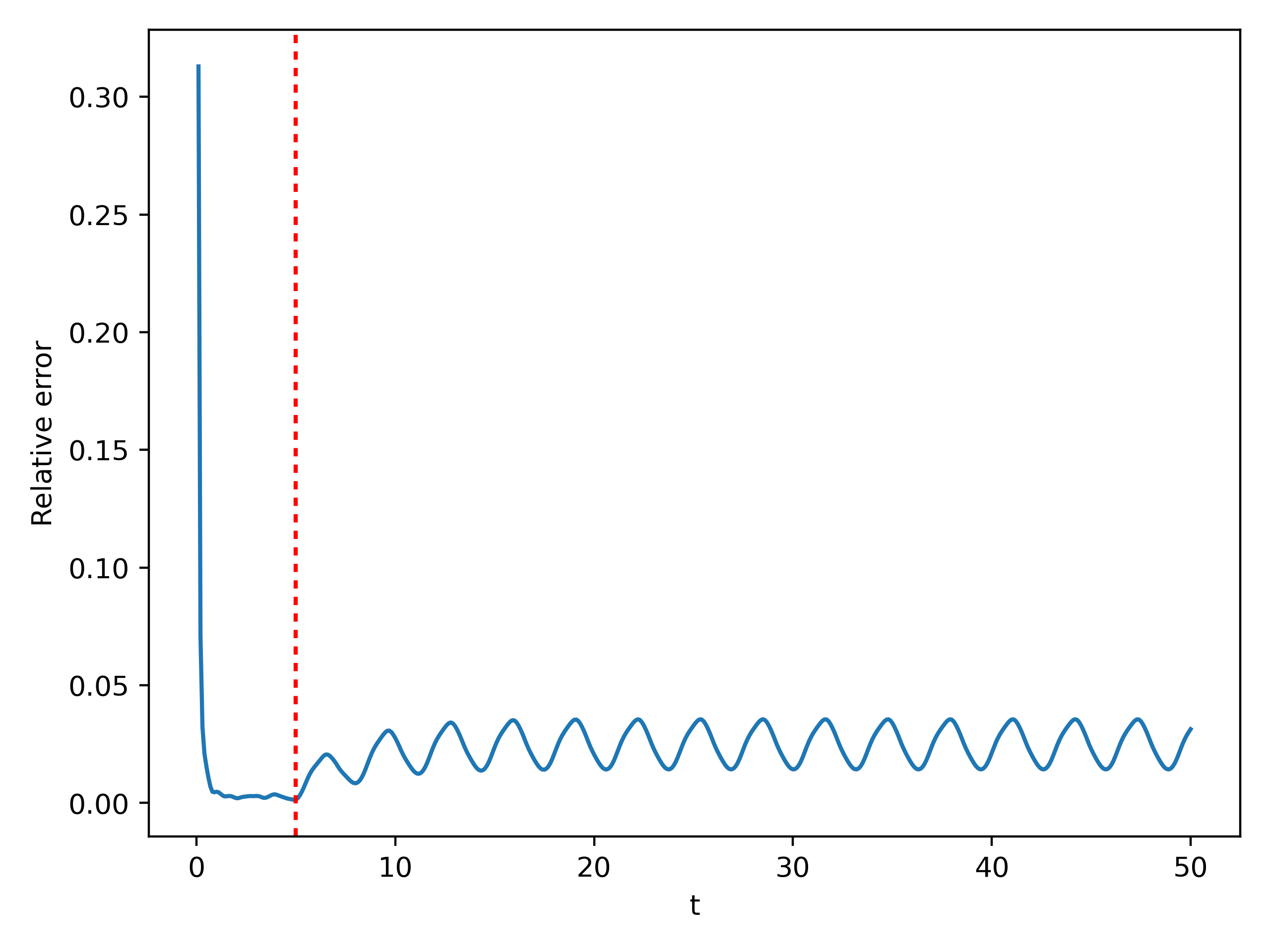}
		\caption{The state error $e_{\Omega}(t)$ over time. The red line indicates the end of the training interval $(0, T)$.}
		\label{fig:examples:advection-diffusion:state:error-over-time}
	\end{figure}

        To evaluate the temporal extrapolation abilities of the model we
        compute the state error $e_\Omega(t)$ beyond the training interval $(0, T)$
        with the results shown in \cref{fig:examples:advection-diffusion:state:error-over-time}.
	We note that the relative error is largest close to $t = 0$ as the state is zero here.
	Beyond the training interval, for $t > T$ we observe that the error stays below 5\% and oscillates with the same period as the inlet $\sin^2(t)$.
	Note that even though we only plot to $t = 50$, we observe no significant deviations from this pattern even at $t = 1000$.

\begin{figure}
		\centering
		\centerline{\includegraphics[width=0.9\linewidth]{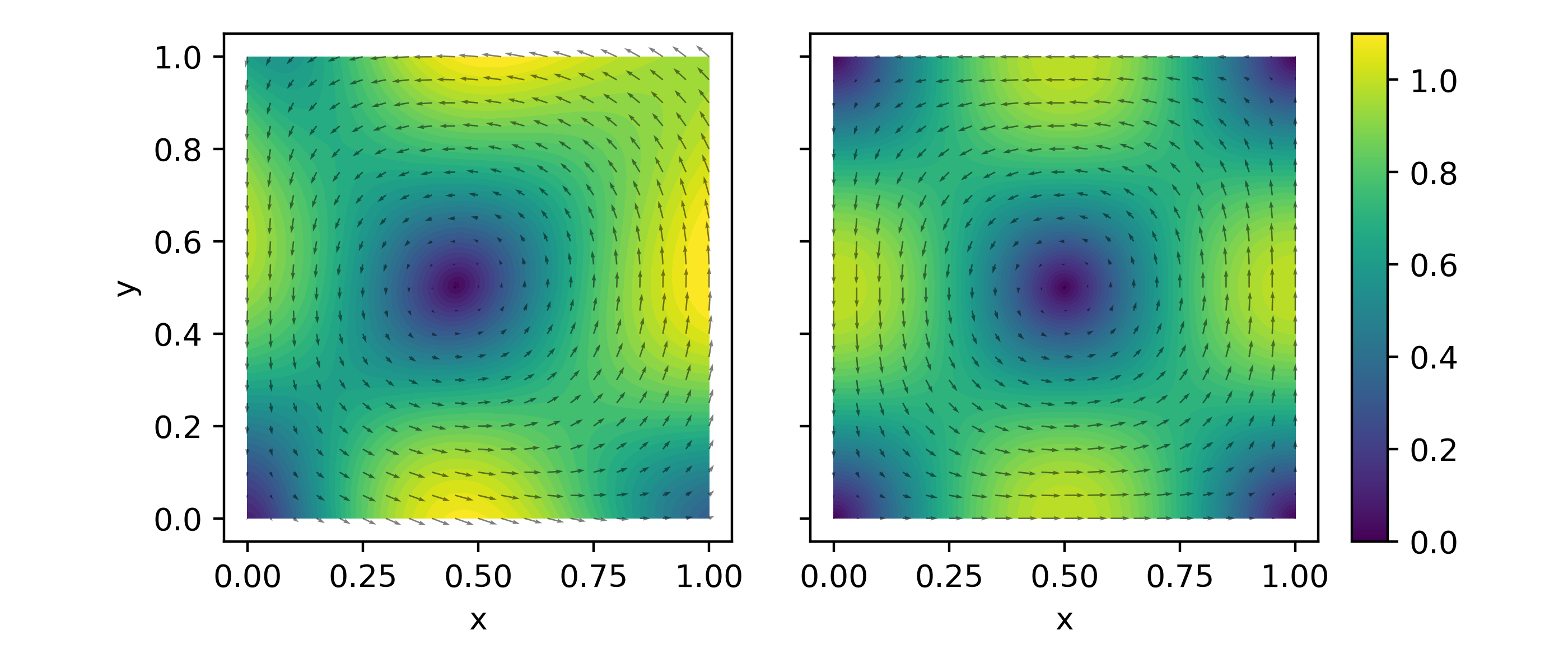}}
		\caption{Left: the neural network prediction over the spatial domain $\Omega$ with constant $u = 0$ to extract the predicted velocity field.
			Right: ground truth advective velocity field over the spatial domain $\Omega$.}
		\label{fig:examples:advection-diffusion:velocity-plot}
	\end{figure}

	We can get a picture of how well the neural network approximates the velocity $\omega$ in
        the advection term $\nabla\cdot(\omega u)$ by supplying artificial state gradients $\nabla u$.
	Using $u = 0$, $\nabla u_1 = (1, 0)$, and $\nabla u_2 = (0, 1)$ \cref{fig:examples:advection-diffusion:velocity-plot} shows the extracted velocity $[\mathcal{N}(x, u, \nabla u_1; W), \mathcal{N}(x, u, \nabla u_2, x; W)]$ along the true 
	velocity $\omega$.
	The predicted velocity field resembles the true one,
	however, there are some distinct differences.
	For one, as opposed to the true velocity, the divergence of the predicted velocity is not zero.
	Further, the predicted velocity field is not tangential to all boundaries, only the inlet boundary.
	
	Note that the artificially supplied arguments used to construct $\hat \omega$ can potentially be tuned to be closer in line with the training data in order to improve accuracy of the predicted velocity.
	However, we did not do this here as in order to assess the better approximation one would need the ground truth, which we assume is unavailable.
	
\begin{figure}
		\centering
		\centerline{
		\includegraphics[width=0.6855\linewidth]{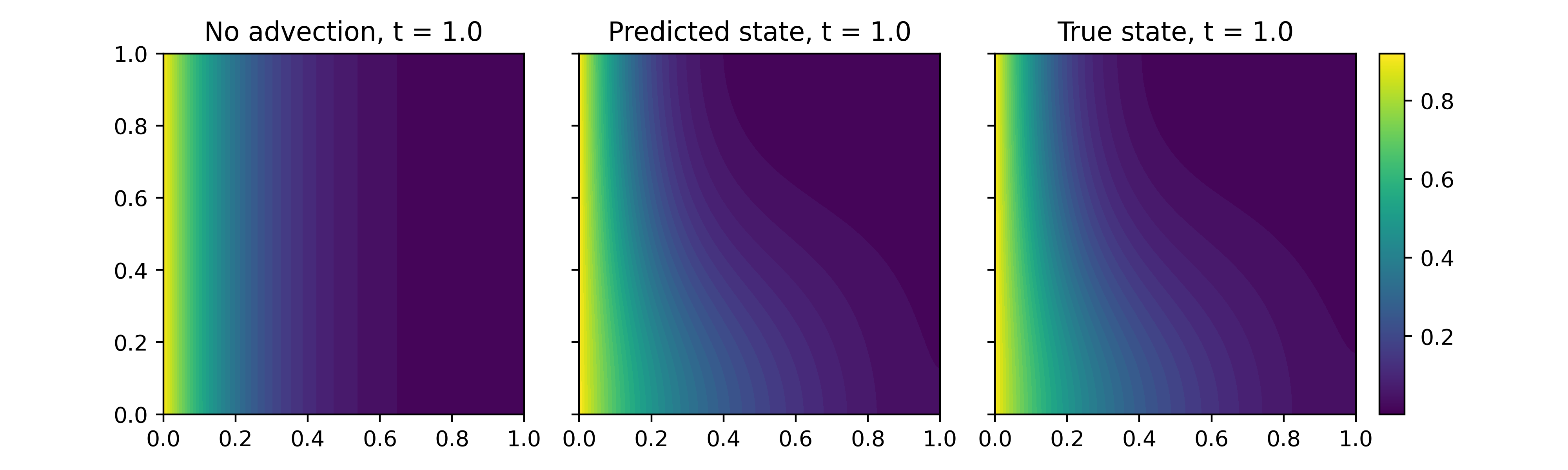}}
		\centerline{
		\includegraphics[width=0.6855\linewidth]{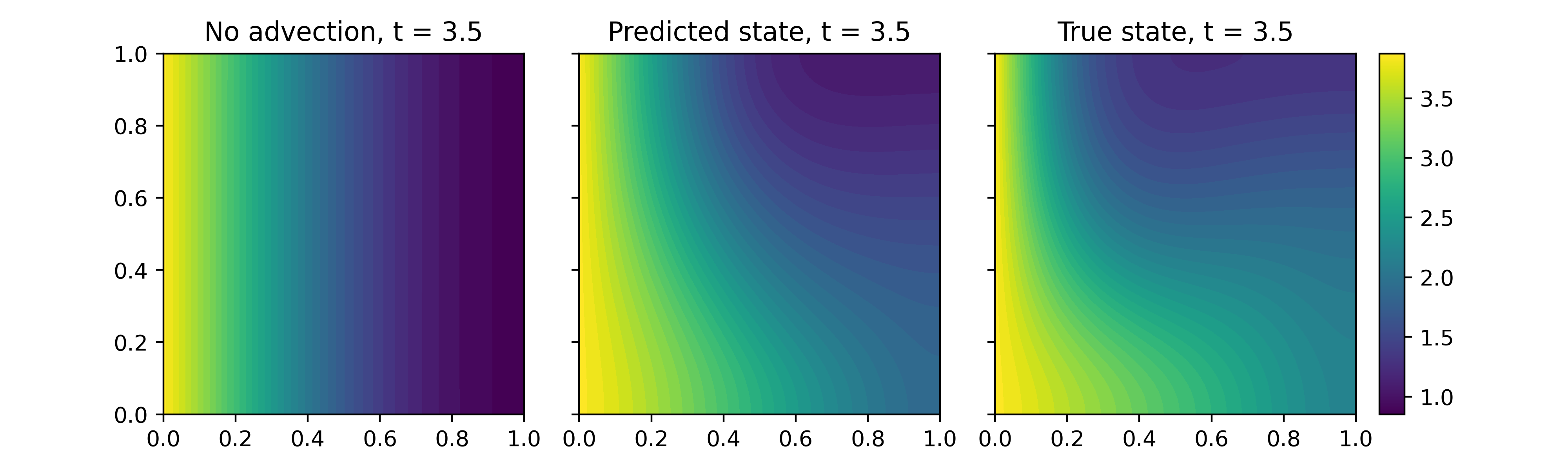}}

		\caption{The predicted and true solution of \cref{eq:examples:advection-diffusion:approx} with a new (unseen) Dirichlet boundary condition with amplitude $2$. On the left the solution using only the PDE part is included ($\mathcal{N} \equiv 0$).
		  Top: the predicted and true states at $t = 1.0$.
		  The relative error, $e_{\Omega}(1.0)$, is 0.44\%.
          Bottom: the predicted and true states at $t = 3.5$.
		  The relative error, $e_{\Omega}(3.5)$, is 11.8\%.}
		\label{fig:examples:advection-diffusion:state:test0:1.0+3.5}
	\end{figure}

	To further investigate the accuracy of the trained network, we consider the same equation with a new Dirichlet boundary condition on the inlet boundary.
	The amplitude is increased from 1 to 2 and the frequency reduced to arrive at $2 \sin^2(0.5t)$.
	While the reduced frequency reduces the absolute value of the gradient, it increases the total quantity of substance coming through the inlet over the time interval $(0, T)$.
	Using the neural network that was trained on the old Dirichlet boundary condition, we solve \cref{eq:examples:advection-diffusion:approx} and measure the error against the ground truth.
	In \cref{fig:examples:advection-diffusion:state:test0:1.0+3.5} the predicted solutions at $t = 1.0, 3.5$ are illustrated.
	The relative error on the test set is $e_{\Omega \times \mathcal{T}} = 6.19\%$ for the full model, and $\mathcal{E}_{\Omega \times \mathcal{T}} = 32.0\%$ for the sub-physics.

    \subsection{Cardiac electrophysiology}
    \label{sec:examples:cardiac}
    The previous examples  used simple fully-connected neural networks 
    and relied on FEniCS/UFL for their implementation. For greater flexibility in defining more complicated NNs it is desirable to take advantage of a dedicated framework. In this final example we  demonstrate the use of more complex neural networks in a hybrid FEM-NN setting where the NNs are implemented in PyTorch \cite{NEURIPS2019_9015}.    
    
    Following \cite[Ch 2.2]{sundnes2007computing} we consider a mathematical model of
    excitable cardiac tissue describing evolution of the transmembrane potential
    $v$ by the coupled PDE-ODE system
    \begin{subequations}\label{eq:examples:cardiac:problem}
      \begin{equation}          
	v_t - \nabla \cdot (M \nabla v) = I_s - I_{\text{ion}}(v, s) \quad \text{in } \Omega \times (0, T),\label{eq:monodomain}
      \end{equation}
      \begin{equation}          
	  s_t = F(v, s) \quad \text{in } \Omega \times (0, T).\label{eq:cardiac_ode}
      \end{equation}
    \end{subequations}
    Here $M$ is the conductivity tensor of the monodomain equation \cref{eq:monodomain} and $I_s$
    is a given stimulus current. The term $I_\text{ion}$ is the ionic current describing a response
    of the cardiac cell with the (internal) state $s(x, t)$ which evolves according to \cref{eq:cardiac_ode}.
    For this example, we let $\Omega = [0, 25]^2$ and $T = 30$.
    We let the conductivity tensor $M$ be the same as the intracellular conductivity ($M_i$) in \cite[Sec 6.2]{farrell_automated_2017}: $M = \mathrm{diag}(g_\mathrm{if}, g_\mathrm{if})$, with $g_\mathrm{if} = 0.174/(\chi C_m)$, $\chi = 140$, and $C_m = 0.01$.

    The functions $I_\text{ion}$ and $F$ are specific to a given cell model
    \cite{farrell_automated_2017} as they describe the cell's dynamics. 
    However, which properties of the cell are relevant for its electrical
    behavior and how those properties interact, is a difficult problem and there
    exists a plethora of cardiac models, see eg. \cite[Ch 2.4]{sundnes2007computing}.
    It is thus reasonable to consider a setting where $F$ and $I_\text{ion}$
    would be learnt from data of a specific patient. Therefore, in the following $F$ and $I_\text{ion}$ are assumed unknown.
    However, we shall maintain the coupling between the terms as in \cref{eq:examples:cardiac:problem}.
    That is, $F$ and $I_\text{ion}$ interact via state vector $s$ and are
    tied to the monodomain equation by the transmembrane potential.
    
    For our experiments, the model parameters, observations and ground truth are determined by using the Hodgkin-Huxley model \cite{hodgkin1952quantitative}.
    The Hodgkin-Huxley model describes a large squid nerve cell, and is therefore not directly applicable to heart models.
    Thus, the action potentials seen in this section are not physiologically meaningful.
    However, the Hodgkin-Huxley model is representative for how cell models are constructed and was later adapted to describe Purkinje fibre cells in the heart \cite[Ch 2.4]{sundnes2007computing}\cite{noble1962modification}.
    
    The decoupled problem is considered, in which the PDE and ODEs are solved separately through a Marchuk-Yanenko splitting scheme.
    Specifically, in order to advance from $t_n$ to $t_{n+1}$ with the timestep $\Delta t = t_{n+1} - t_{n}$, the PDE 
	\begin{align*}
		\frac{v^{n+1} - v^n}{\Delta t} - \nabla \cdot (M \nabla v^{n+1}) &= I_s - I_{\text{ion}}(v^n, s^n)\quad\mbox{ in }\Omega,
	\end{align*}
	where $v^n = v(\cdot, t_n)$ and $s^n = s(\cdot, t_n)$,
        is first solved using finite elements.
	Then, the system of ODEs is solved using an explicit Euler scheme
	\begin{align*}
		s^{n+1} = s^n + \Delta t F_i(v^{n+1}, s^n).
	\end{align*}
	Synthetic data was generated using this scheme on finite elements and a uniform timestep of $\Delta t = 10^{-3}$.
	
	The goal of this example is to approximate both $I_\text{ion}$ and $F$ using NNs
    trained on observations of only $\{v(x, \tau_i)\}_i$ for a given initial potential $v(x, 0)$ and stimulus function $I_s(x, t)$.
	We construct the hybrid model such that it has the same number of cell states as in the Hodgkin-Huxley model.
    Both NNs take in only $v$ and $s$ evaluated at a point $(x, t) \in \Omega \times (0, T)$,
    and are not provided any spatial or temporal coordinates, nor any derivative information.
	The same scheme as was used to generate synthetic data was used for training the hybrid model.
    However, the time step in training was larger at $\Delta t = 0.01$.
    Furthermore, both neural networks are implemented in PyTorch and interpolated at the finite element mesh vertices.	
	The forcing and ionic current terms are approximated by two feedforward NNs with 3 hidden layers and 16 neurons in each layer. Both NNs use hyperbolic tangent activation functions with a skip connection between the 1st and 3rd layer.
    
	Although estimation of $s(x, 0)$ from observations of $v$ might be possible using a NN as suggested in \cite{ayed_learning_2019}, the hybrid model assumes that $s(x, 0)$ is zero for simplicity.
	Thus, one cannot expect a reconstruction of any meaningful cell states.
	
	The neural networks were trained in two steps, first they were trained using observations at each $\tau_i$ as the input to predict the next state state $\hat v$ at $t = \tau_i + \Delta t$.
	This process is similar to training on subintervals $[\tau_{i-1}, \tau_{i}]$, but allows gradient contributions for the hidden cell state $s$ to flow between the subintervals.
	Second, the neural networks were afterwards further trained while using only the observation at $t = 0$ as input.
	In both training steps we use the mean squared error (MSE) of the transmembrane potential as the loss function, and L-BFGS as the optimization algorithm.
	
	For the training data we apply a stimulus for $3 \leq t \leq 4$ in the upper left corner of the domain:
	\begin{align}
	    I_s(x, t) = \begin{cases}
	        100 (1 - \frac{3}{5}d(x)), \quad \text{for } d(x) < \frac{5}{3}, \, 3 \leq t \leq 4, \\
	        0, \quad \text{otherwise},
	    \end{cases}
	\end{align}
	where $d(x) = d(x_0, x_1) = x_0 + (25 - x_1)$.
	
\begin{figure}
		\centerline{
		\includegraphics[width=0.5\linewidth]{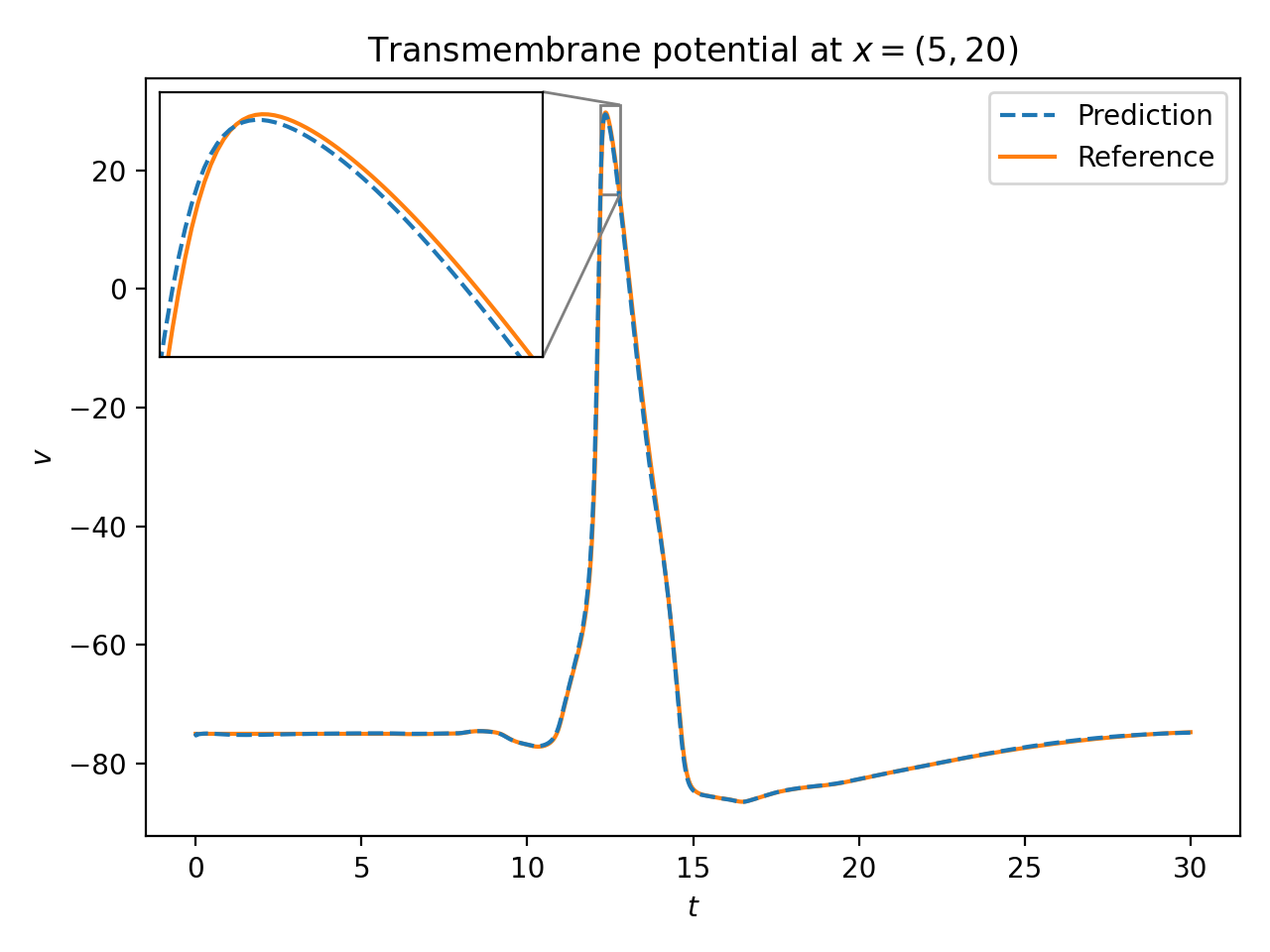}
		\includegraphics[width=0.5\linewidth]{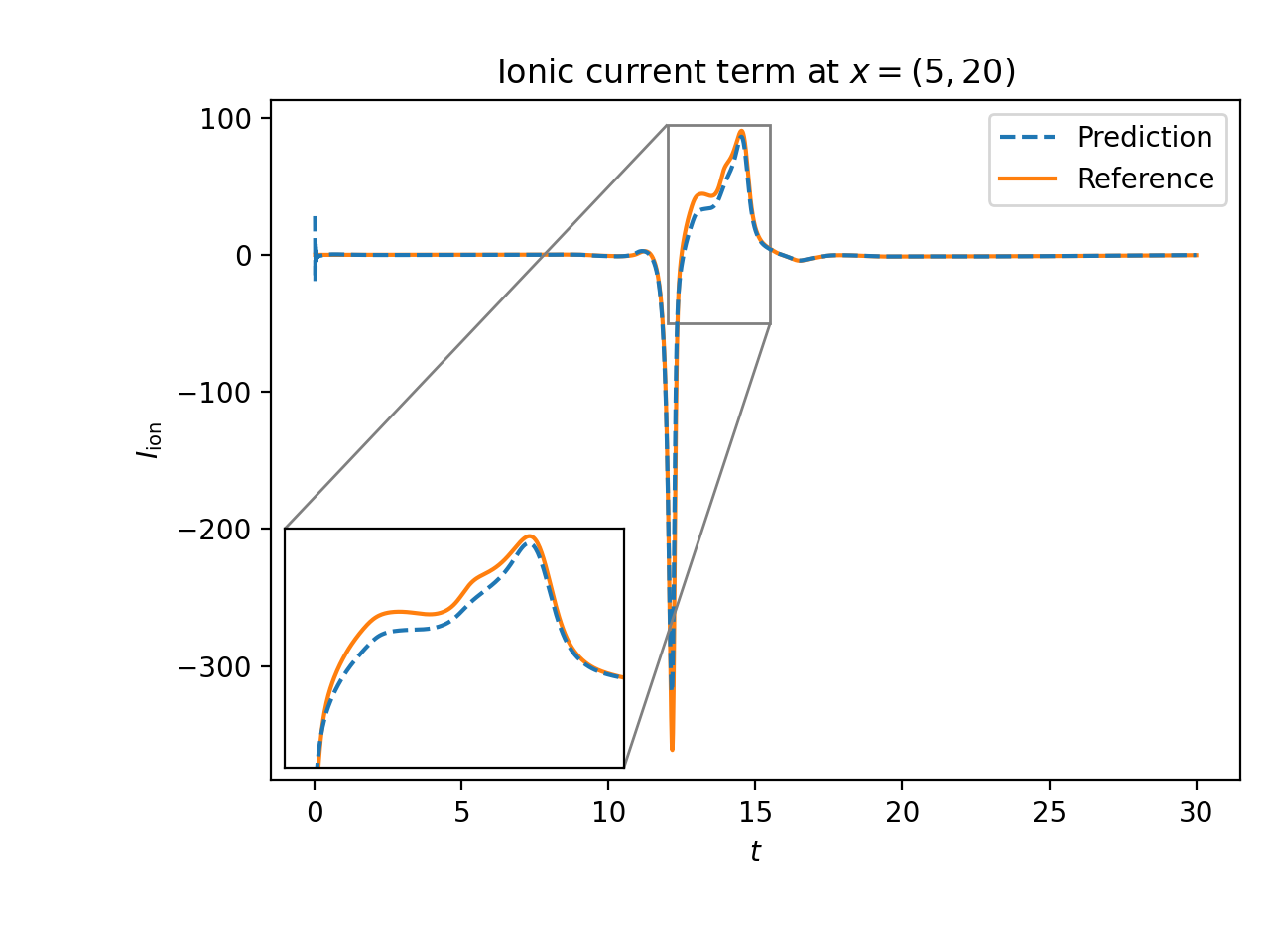}}
		\caption{Left: the predicted and true state $v$ at $x = (5, 20)$ over time.
			Right: the predicted and true $I_{\text{ion}}$ term at $x = (5, 20)$ over time.
			The external stimulus and initial condition are the same as during training.}
		\label{fig:examples:cardiac:training:result:v}
	\end{figure}

\begin{figure}
		\centerline{
		\includegraphics[width=0.9\linewidth]{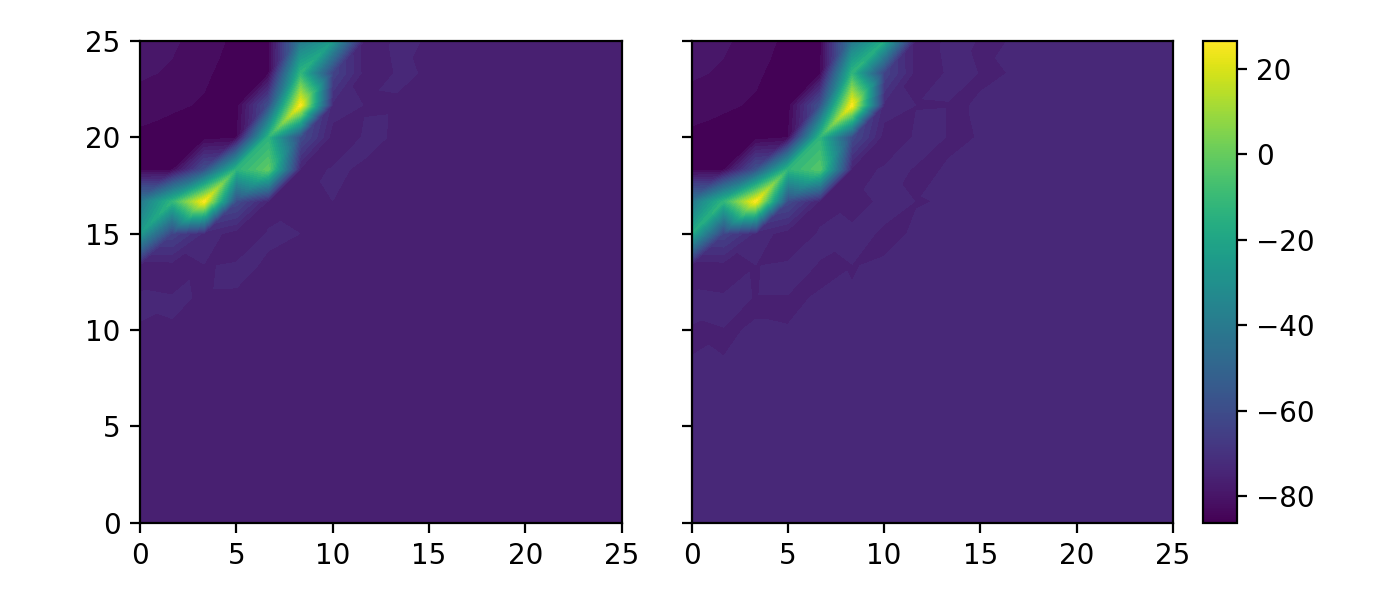}}
		\caption{Left: the predicted transmembrane potential $v$ at time $t = 15.0$.
			Right: the true transmembrane potential $v$ at time $t = 15.0$
			The external stimulus and initial condition are the same as during training.}
		\label{fig:examples:cardiac:training:result:v:omega}
	\end{figure}

	The predicted state after training the NNs can be seen in \cref{fig:examples:cardiac:training:result:v} together with the predicted ionic current term $I_\text{ion}$ at spatial coordinate $x = (5, 20)$.
	The predictions seen in \cref{fig:examples:cardiac:training:result:v} use the same initial condition and external stimulus as during training.
	We observe that both the transmembrane potential $v$ and the ionic current response are fairly accurately predicting the reference.
	We make the same observation when looking at the spatial distribution of the transmembrane potential at time $t = 15.0$ in \cref{fig:examples:cardiac:training:result:v:omega}.
	Quantified, the prediction error in the whole domain for the transmembrane potential is $e_{\Omega \times (0, T)} = 0.48\%$, however, this is misleading due to the equilibrium of the state $v$ not being at 0.
	The absolute error of the predicted transmembrane potential over the whole spatio-temporal domain is $441.5$.
        
\begin{figure}
		\centerline{
			\includegraphics[width=0.5\linewidth]{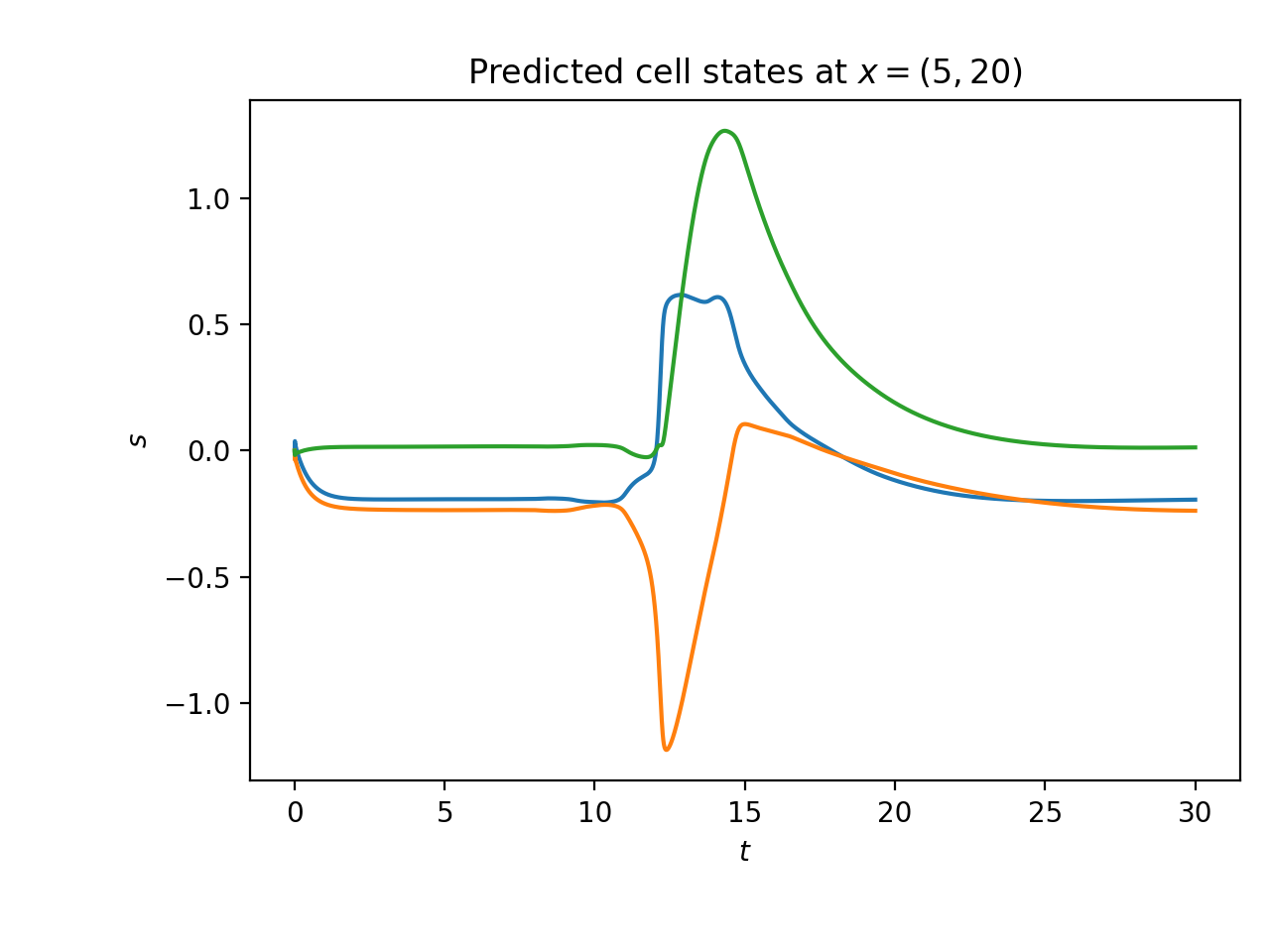}
			\includegraphics[width=0.5\linewidth]{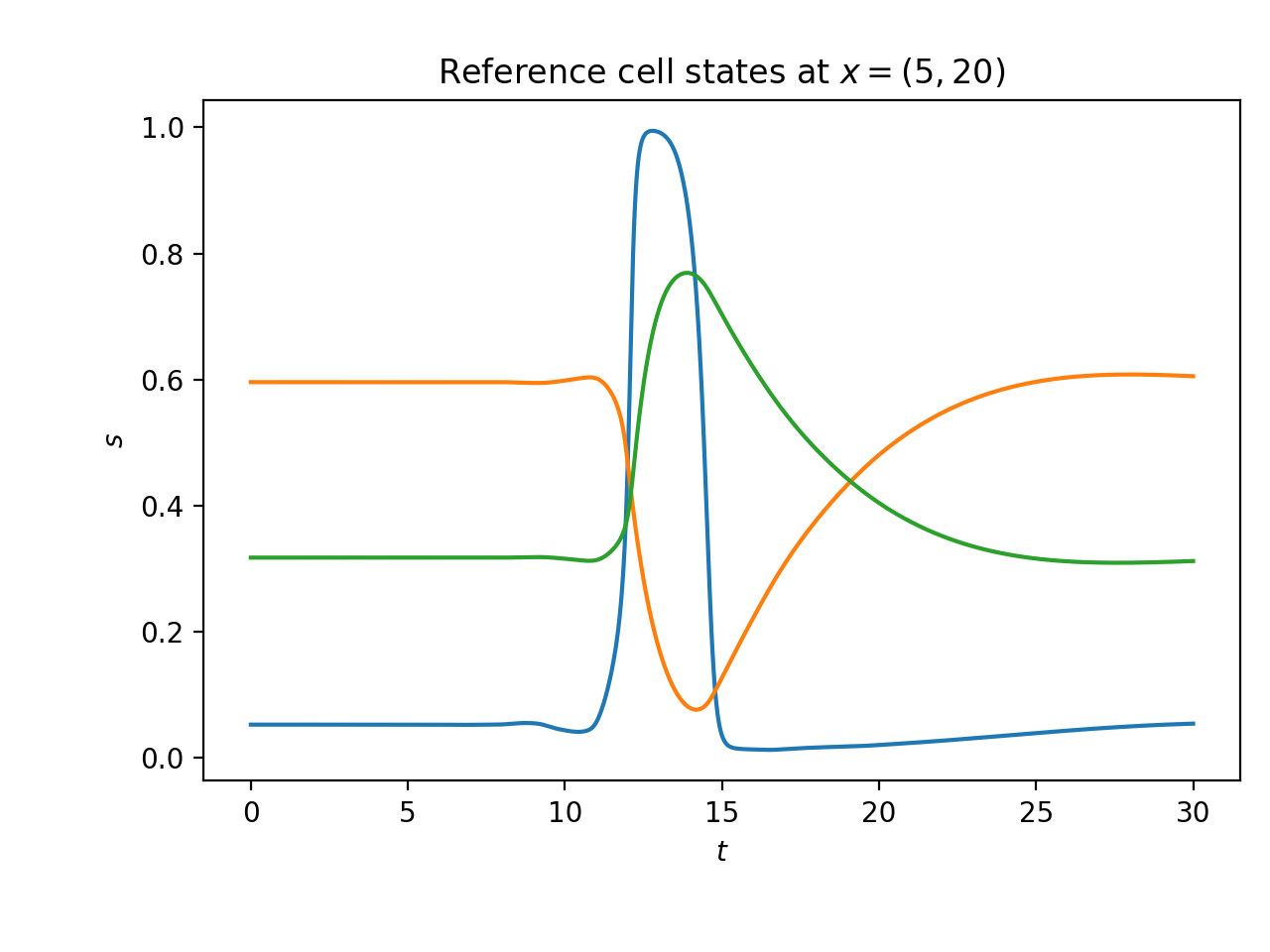}}
		\caption{The hidden cell states $s$ in the trained hybrid FEM-NN model (left) and the Hodgkin-Huxley model (right) at the spatial point $x = (5, 20)$.}
		\label{fig:examples:cardiac:training:result:s}
	\end{figure}

	We also inspect the predicted hidden cell states for $s$, and compare to the true $s$ from the Hodgkin-Huxley model. The values of these hidden states at spatial coordinate $x = (5, 20)$ are plotted in \cref{fig:examples:cardiac:training:result:s}.
	The predicted cell states are all initialized to zero, but quickly transitions to an equilibrium until the wave triggered by the external stimulus reaches $x$.
	While the scales of the approximate and true model are quite different, the predicted cell states resemble the true cell states.
	This is surprising given that the NN has not been exposed to the true cell states during training and is not explicitly constrained to any particular dynamics for the cell state.
	 
\begin{figure}
		\centerline{
			\includegraphics[width=0.5\linewidth]{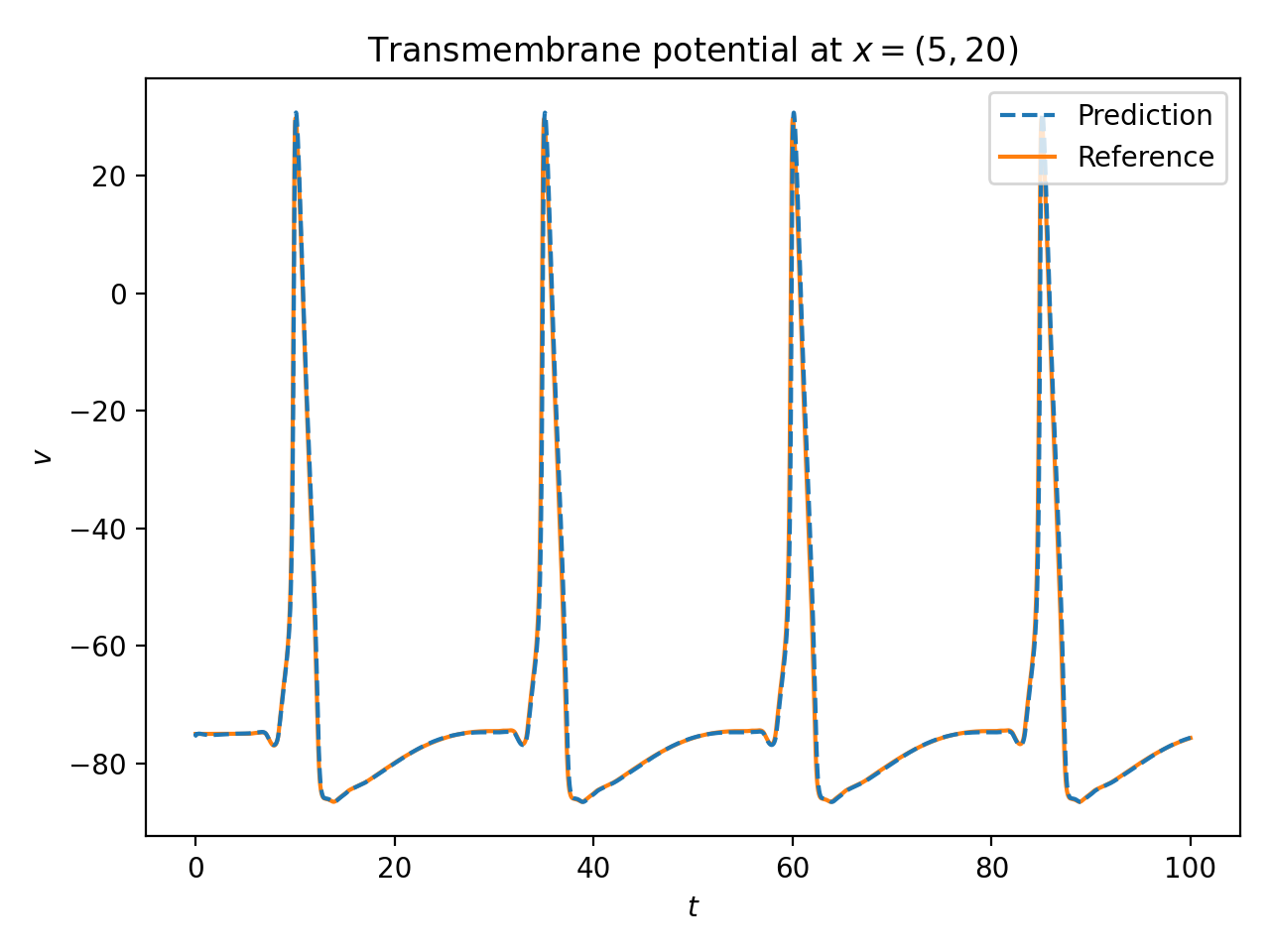}}
		\caption{The predicted and true state $v$ at $x = (5, 20)$ over time when using a new (unseen) stimulus applied in a box at the center of $\Omega$ repeating every 25 units of time.}
		\label{fig:examples:cardiac:test:result:v}
	\end{figure}

\begin{figure}
		\centerline{
			\includegraphics[width=0.9\linewidth]{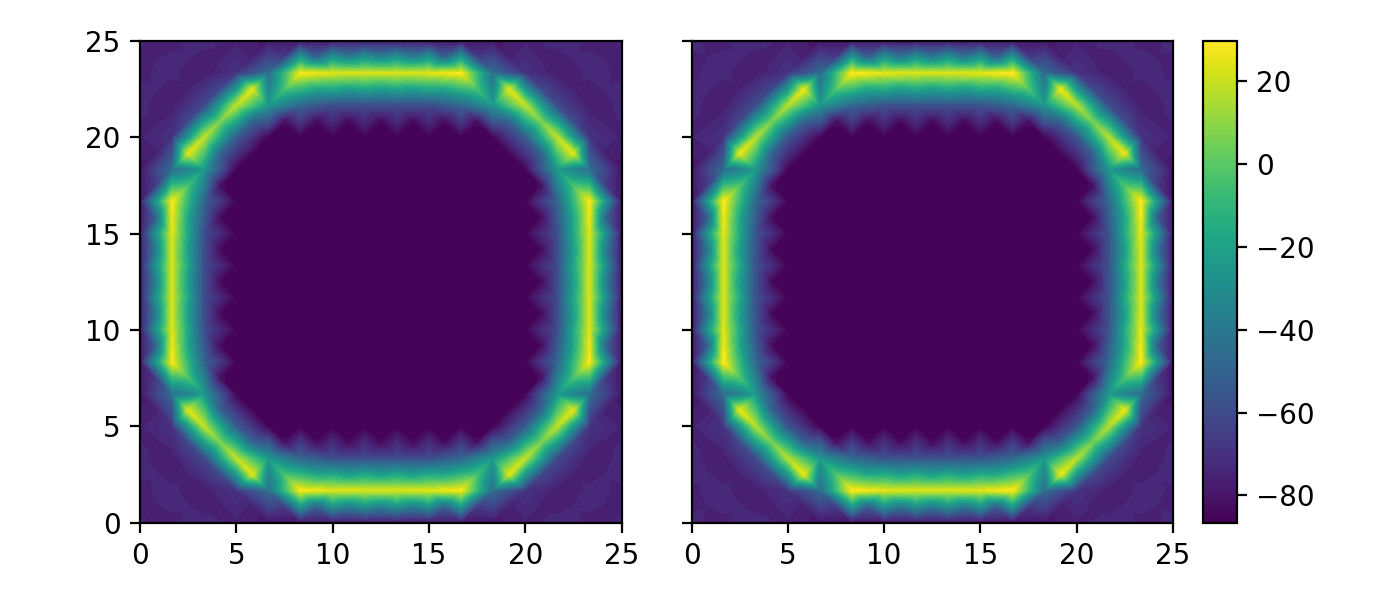}}
		\caption{Left: the predicted transmembrane potential    $v$ at time $t = 12.0$.
			Right: the true transmembrane potential $v$ at time $t = 12.0$.
			The external stimulus was not seen during training, and applies in a square in the center of the spatial domain $\Omega$.}
		\label{fig:examples:cardiac:test:result:v:omega}
	\end{figure}

	The trained neural networks were tested on a new stimulus applied at the center of the domain $\Omega$ and at a different time $t \in [5.0, 6.0]$ after which it repeats every 30 units of time.
	This stimulus is applied as a constant $I_s = 100$ in a centered square with length $10$.
	The result of the predictions of the trained model at the spatial point $x = (5, 20)$ plotted over time can be seen in \cref{fig:examples:cardiac:test:result:v}.
	The trained networks provide accurate predictions of $v$ also for this new stimulus.
	Thus, the trained FEM-NN model is able to extrapolate in time for this stimulus, generalizing beyond a single action potential that was seen during training.
	We also look at the predictions over all of $\Omega$ at a snapshot $t = 12.0$ in \cref{fig:examples:cardiac:test:result:v:omega}, and observe that the predictions are accurate throughout $\Omega$.
	The relative error on the test set is $e_{\Omega \times (0, 100)} = 1.7\%$ with an absolute error of $2462$.
	
    \section{Conclusion}\label{sec:conclusion}
    
    This paper discussed hybrid FEM-NN
    models, in which PDEs are augmented with neural networks to represent unknown terms, and discretised using the finite element method.
    The methodology was demonstrated on several examples 
    of inverse problems and also compared to other approaches.
    When solving the inverse heat conduction problem with a discontinuous 
    diffusion coefficient we observed that the hybrid method produces more 
    accurate approximations of the coefficient (with less overshoots at the point of the jump 
    in the material) compared to both the pointwise estimation and the standard PINNs. 
    Compared to PINNs fewer optimization steps were required to reach this accuracy 
    suggesting that the strong enforcement of the PDE in the hybrid FEM-NN model improves the speed of convergence.
        We also demonstrated the recovery of an unknown advection term $u\mapsto \omega\cdot \nabla u$ in the advection-diffusion equation. The recovered 
        operator was analyzed in terms of the extracted approximated advection velocity. We found that the velocity matches reasonably well with the ground truth.
    Lastly, we demonstrated the use of hybrid FEM-NN to reconstruct a cell model in the monodomain equation.
    The trained model matched well with the Hodgkin-Huxley model used for data, and could generalize to different spatially varying stimuli.

    While the hybrid FEM-NN approach offers opportunities and flexibility, it also introduces new challenges: (i) the 
    nonlinear and linear subproblems arising from the FEM-NN discretisation must remain solvable throughout training 
    and (ii) well-established PDE solver tools are not optimised for large batch training or repeated evaluations on GPU, and (iii) the variational formulation used in the FEM requires spatial integration over the NN, for which efficient quadrature rules are currently unknown.
    
    The learned advection term is relatively simple given that the gradient of the solution was supplied. An interesting avenue for future work is therefore to apply the hybrid FEM-NN models to learn advection or other physical terms for highly nonlinear problems.

    Furthermore, all neural networks shown have been relatively small, and finding optimal designs of the neural networks have not been part of our work.
    While all the tested NN are global functions mapping each point separately, avoiding the discretisation to be designed into the network, this also comes with a limitation with regards to learning differential operators.
    Because the network will only have local information for its input functions, it cannot produce differential terms without derivatives being explicitly supplied.
    Therefore, investigating the capabilities of a local NN defined on stencils, possibly designed as a convolutional NN as in \cite{rackauckas_universal_2020}, is an interesting pathway forward.
    
    Finally, we have observed that the trained NN will attempt to compensate and produce 
    a more accurate estimate without discretisation when during the training the discretisation 
    error is high (the results are not included in the paper). In addition, even if the discretisation 
    scheme is unstable in the ground truth model, the NN seems to be able to produce outputs that give 
    reasonable solutions. The observation poses an interesting question for future work, namely: Is it 
    possible to train a NN to alleviate discretisation errors and instabilities in such a way that it generalizes?
    
    \section{Acknowledgments}\label{sec:acknowledgments}
    This work was supported by the Research Council of Norway through the FRINATEK program, project numbers 303362.
    Sebastian Mitusch was supported by the Norwegian Ministry of Education and Research.
    Miroslav Kuchta acknowledges support from the Research Council of Norway grant no. 280709.
    
    \bibliography{paper}

\end{document}